\documentclass[12pt]{amsart}
\usepackage{amsmath}
\usepackage{amssymb}
\usepackage{amsfonts}
\usepackage{amsthm}
\usepackage{verbatim}
\usepackage{amscd}
\usepackage{cite}
\usepackage{leftidx}
\usepackage{enumerate}
\usepackage{txfonts}
\usepackage{manfnt}
\usepackage{amscd}
\usepackage[mathscr]{eucal}
\usepackage{hyperref}
\usepackage{datetime2}
\usepackage{graphics}
\usepackage{xcolor}
\usepackage{bookman}
\def\cZ{\mathcal Z}

\newtheorem{thm}{Theorem} 

\newtheorem*{thm*}{Theorem}
\newtheorem*{prop*}{Proposition}
\newtheorem{cor}[thm]{Corollary}
\newtheorem*{cor*}{Corollary}

\newtheorem{lem}[thm]{Lemma}
\newtheorem*{lem*}{Lemma}

\newtheorem*{claim*}{Claim}
\newtheorem{prop}[thm]{Proposition}

\newtheorem{defn}[thm]{Definition}
\theoremstyle{remark}

\newtheorem{rem}[thm]{Remark}
\newtheorem*{rem*}{Remark}
\newtheorem{crit-rem}[thm]{Critical remark}
\newtheorem{remarks}[thm]{Remarks}

\newtheorem{example}[thm]{Example}
\newtheorem*{example*}{Example}
  
\newtheorem*{defn*}{Definition}
\newtheorem*{con*}{Conjecture}

 \DeclareMathOperator{\diag}{Diag}

 \DeclareMathOperator{\gr}{Gr}
\def\inv{^{-1}}

\DeclareMathOperator{\Ext}{Ext} \DeclareMathOperator{\Hom}{Hom}
\DeclareMathOperator{\coker}{coker}

\def\cC{\mathcal C}

\def\refp #1.{(\ref{#1})}

\newcommand\ceil [1] {\lceil #1 \rceil}

\newcommand{\A}{\mathcal{A}}

\newcommand{\remainder}[2]{{#1}\text{\%}{#2}}

\newcommand{\emphred}[1]{\emph{\color{red}#1}}

\def\sbr #1.{^{[#1]}}
\def\sfl #1.{^{\lfloor #1\rfloor}}

\def\inv{^{-1}}
\def\?{{\bf{??}}}

\newcommand{\PGL}{\rm PGL}

\def\A{\mathbb A}

\def\C{\mathbb C}
\def\P{\mathbb P}

\def\Z{\mathbb Z}

\def\O{\mathcal O}

\def\rk{\text{rk}}

\def\g{\mathfrak g}

\def\gl{\mathfrak{gl}}

\def\m{\mathfrak m}
\def\k{\mathfrak k}

\def\hom{\mathfrak {hom}}

\def\1/2{\frac{1}{2}}

\def\simto{\stackrel{\sim}{\rightarrow}}
\def\2{{[2]}}
\def\l{\ell}
\def\nl{\newline}

\def\hom{\mathcal{H}\mathit{om}}

\def\<{\langle}
\def\>{\rangle}

\def\2{{[2]}}
\def\l{\ell}

\def\scl #1.{^{\lceil#1\rceil}}
\def\spr #1.{^{(#1)}}
\def\sbc #1.{^{\{#1\}}}

\def\subpr#1.{_{(#1)}}

\def\beq{\begin{equation*}}
\def\eeq{\end{equation*}}

\newcommand{\mlog}[1]{\langle-\log {#1} \rangle}
\def\g3{{\Gamma\spr 3.}}

\newcommand{\eqspl}[2]{
\begin{equation}\label{#1}
\begin{split}
#2\end{split}\end{equation}}

\newcommand{\exseq}[3]{
0\to #1\to #2\to #3\to 0
}

\newcommand{\beginalphaenum}{
\begin{enumerate}\renewcommand{\labelenumi}{ }
\item \begin{enumerate}
}

\def\eex{\end{rm}\end{example}}

\newtheorem{examples}[thm]{Examples}
\begin{document} 
\setcounter{thm}{-1}
	\title{Holonomy and bundle stability\\ for curves
	   on Fano hypersurfaces}
	\author 
	{Ziv Ran}
%
%
	\thanks{arxiv.org 2211.12661} 
	\date {\DTMnow}
%
%
	\address {\nl UC Math Dept. \nl
	Skye Surge Facility, Aberdeen-Inverness Road
	\nl
	Riverside CA 92521 US\nl 
	ziv.ran @  ucr.edu\nl
	\url{https://profiles.ucr.edu/app/home/profile/zivran}
	}
	
	 \subjclass[2010]{14n25, 14j45, 14m22}
	\begin{abstract}
	For every $n\geq 3,  g\geq 2$
	 and all large enough $e$ depending on $n,g$, there exist curves
	of genus $g$, degree $e$ in  a general hypersurface of degree $n$ in $\P^n$, or in $\P^n$ itself,
	 whose normal bundle $N$  is stable, as is any sufficiently general full-rank subsheaf of $N$.
	 For $g=1$, $N$ is semi-stable.
On general hypersurface of degree $d< n$ in $\P^n$, such that
a certain arithmetical condition on $d,n, g $  holds, there exists an arithmetical
progression of $e$ values so that curves of degree $e$ and  genus $g$ with
semistable normal bundle exist.
Previous results were restricted to certain cases
with  ambient space $\P^n$. A main tool in the proof is an analogue of 
parallel transport (resp. holonomy) for
bundles on a chain (resp. cycle)  of rational curves. 
		\end{abstract}
	\maketitle
	\section*{Conflict of Interest Statement}

	This manuscript is not under review at another journal or other publishing venue.
	
	\par	  The author has no affiliation with any organization with a direct or indirect 
	financial interest in the subject matter discussed in the manuscript
	\section*{Data Availability Statement}
	There is no data set associated with this paper
	\section*{AI statement}
	This paper was written solely by the above-named author; no part was played by any AI.

	\setcounter{section}{-1}
	\setcounter{enumi}{-1}
	\section{Introduction}
	\subsection{Set-up}
	A vector bundle $E$ on a curve $C$ is said to be semistable if every subbundle $F\subset E$
	has slope $\mu(F):=\deg(F)/\rk(F)\leq \mu(E)$. It is a natural question whether
	for a sufficiently general curve $C$ on a given projective variety $X$, 
	with given homology class or degree $e$ and genus $g\geq 1$, the normal
	bundle $N=N_{C/X}$ is semistable. This intuitively means $C$ moves evenly in all directions.
	\subsection{Known results}
Except for the case where $\mu(N)$ is an integer, existence results for curves with semistable normal bundle
have been sporadic and  restricted to the case of ambient space $X=\P^n$, mostly $n=3$.
They are due to many 
mathematicians 	since the early 1980s,
	including Ballico-Ellia \cite{ballico-ellia-normal},
	Coskun-Larson-Vogt \cite{coskun-larson-vogt-normal},
	Ein-Lazarsfeld \cite{ein-laz-normal}, Ellia \cite{ellia-normal}, 
	Elligsrud-Hirschowitz \cite{ellingshi}, Ellingsrud-Laksov \cite{ellingslak},
	Ghione-Sacchiero \cite{ghione-sacchiero-normal}, Eisenbud-Van de Ven 
	\cite{eisenbud-vandeven-normal},  Hartshorne \cite{hartshorne-space-curves},
	Hulek \cite{hulek-elliptic}, 
	Newstead \cite{newstead-normal}, Sacchiero \cite{sacchiero-stable-normal}. Notably, Ein and Lazarsfeld \cite{ein-laz-normal} have shown
	that an elliptic curve of (minimal) degree  $n+1$ in $\P^n$ has semistable normal bundle, and
	Coskun-Larson-Vogt \cite{coskun-larson-vogt-normal} have shown 
	recently that with few exceptions,
 Brill-Noether general curves in $\P^3$ have stable normal bundle.
\par	Along different lines, it was shown by Atanasov-Larson-Yang \cite{alyang} that the normal bundle
	to a general nonspecial curve in $\P^n$ is interpolating or balanced in the sense that for a general
	effective divisor $D$ of any degree one has either $H^0(N(-D))=0$ or $H^1(N(-D))=0$. It is known
(see lemma \ref{int-slope-lem})	that for any bundle with \emph{integral} slope, interpolation implies semistability. Thus, \cite{alyang}
	yields curves with semistable normal bundle in $\P^n$ with any genus $g$ and  degree $e$ such that
	$2e+2g-2\equiv 0\mod n-1$.\par 
	For curves on hypersurfaces $X$ of degree $d\leq n$ in $\P^n$, 
the results of \cite{elliptic} on interpolation again yield some curves whose normal bundle is balanced
with integer slope, hence semistable. In the case of anticanonical hypersurfaces ($d=n$) this includes curves on any genus $g$ 
and sufficiently large degree $e$ such that $e+2g-2\equiv 0\mod n-2$ (for $d<n$ there are additional
congruence conditions). Apparently there are no known results where $N$ has non-intergal slope.
\subsection{New results}
	The purpose of this paper is to expand
	the collection of curves  with known stable or semistable normal bundle on hypersurfaces.
To compactify the statement, we will adopt the following terminology. 
A $(g,e)$ curve is one of genus $g$ and degree $e$;  a curve $C$ is \emph{normally semistable}
(resp. \emph{stable}) in a given ambient
space $X$ if its normal bundle $N_{C/X}$ is semistable (resp. stable).
$C$ is \emph{ambient (semi) stable} if the restricted tangent bundle $T_X|_C$ is.
A bundle is said to be hyper-stable if its sufficiently general subsheaf of any given
finite colength is stable; ditto for semi-, normally- etc.
 \par
Our new results are as follows. First for semistability
(see Theorems \ref{genus g}, \ref{ambient-semistable-thm},
\ref{hypersurface-thm},  \ref{lower}
for precise statements).
\begin{thm*}[(Semi)stability] (a) Given $g\geq 1,e>>0, 3\leq n,d\leq n$, then\par
(i) a general $(g,e)$ curve in $\P^n$ is normally hyper-stable (resp. semistable) if $g\geq 2$
(resp. $g=1$) and ambient semi-stable for any genus $g\geq 1$; \par
 (ii) a general hypersurface of degree $n$ in $\P^n$ contains
a $(g,e)$ normally hyper-stable (resp. semi-stable) curve if $g\geq 2$ (resp. $g=1$); \par
(iii) if $X$ be a general hypersurface of  degree $d< n$ 
in $\P^n$, and $g\geq 1$ is such that the divisibility condition 
\[	{((d-2)(n-d+1), d(n-2)\ \ |\ \ 2(n-d)(g-1).
			}\]
			 holds.
Then there exists an arithmetic progression of $e$ values so that  $X$ contains
a normally semistable curve of degree $e$. 
For example, the divisibility  condition holds for  $g= 1$ and also,
 \par  - for  $d=n-1$, if either $n$ is even or  $g$ is odd; 
\par
- for $d=n-2$, if either  $d\not\equiv 2\mod 3$ or $g\equiv 1\mod 3$;\par
- for all $d, n$, if $g-1$ is sufficiently divisible.
\end{thm*}
\noindent\textbf{Remarks:}\par
\begin{enumerate}\item
The normal bundles in question have in general non-integral slope.	\item
By Riemann-Roch, a $(g, e)$ curve  $C\subset X_d\subset\P^n$ with semistable normal
bundle $N=N_{C/X}$ has $H^1(N)=0$ provided $e>\frac{2(n-3)(g-1)}{n+1-d}$.
Therefore the curves in question are unobstructed.\item
To my knowledge these are the first cases, beyond $g=1$ or $n=3$ or the case
of integral slope, of genus-$g$ curves in $\P^n$ with known semistable normal bundle.\item
Other than the integral slope case, there are no such results
in the literature for hypersurfaces other than $\P^n$.\item
For $n=3$ Coskun-Larson-Vogt have proven stability of the normal bundle for a larger
set of degrees and genera.\item
We do not state a result on stability for curves on hypersurfaces of degree $<n$ in
$\P^n$.\item
There is no reason to expect the range of curve degrees that we obtain to be optimal:
e.g. in the case of $\P^n$ one might optimistically expect (semi)stability for all
Brill-Noether general curves.
\end{enumerate}\qed\par
As for the geometric meaning of a curve
having semistable normal bundle, note that a smooth subvariety
containing the curve yields a subbundle of the normal bundle, so
from the standard formula for the degree of the normal bundle and the definition of semistability,  we have:
\begin{lem*}\label{semistable-lem}
Let $C\subset X$ be a smooth curve of genus $g$ with semistable normal bundle
on a smooth $n$-dimensional variety. 
Then for any subvariety $Y\subset X$
of dimension $m$ containing $C$ and smooth along it, we have
\[(n-1)C.(-K_Y)+(n-m)(2g-2)\leq (m-1)C.(-K_X).\]
\end{lem*}
\subsection{Method of proof} Our strategy is to first prove,  in Theorem \ref{g=1-thm},
semistability for genus 1 using a 'fish fang' degeneration (compare \cite{caudatenormal}, \cite{elliptic})
where the embedding $C\subset\P^n$ degenerates to
\[C_1\cup_{p,q}C_2\subset P_1\cup _EP_2\]
where $P_1, P_2$ are blowups of $\P^n$ in a suitable $\P^m$ resp $\P^{n-1-m}$ with common
divisor $E=\P^m\times\P^{n-1-m}$ and where
$C_1, C_2$ are rational with $C_1\cap E=C_2\cap E=\{p,q\}$. 
In fact we will prove a more general result showing semistability of general modifications
of the normal bundle.
A critical role in the proof is played by a generality property for a
'parallel transport' type map between the fibres of the normal bundle at $p$ and $q$,
and a resulting 'holonomy' at $p$ (going to $q$ along $C_1$ and back along $C_2$).
The holonomy being general enough forces the normal bundle to decompose as a direct
sum of general line bundles.
This result is then used to prove 
the main result for $\P^n$, Theorem
\ref{genus g}, by an induction on the genus, using another fang degeneration.
We will also prove a semistability result for the restriction on the ambient tangent
bundle of $\P^n$ on curves of genus $>0$ and sufficiently high degree. \par
The result for anticanonical hypersurfaces, Theorem \ref{hypersurface-thm} (slightly more general than the above
statement), is proven using the $\P^n$ case
plus a fan-quasi cone degeneration similar to the one used in \cite{caudatenormal}. 
Finally the result for lower-degree hypersurfaces,  Theorem \ref{lower}, is proven using a suitable fang 
as in \cite{caudatenormal}, essentially using  a suitable projective bundle
over a $\P^m$,  using semistablity of the horizontal part of the normal bundle by Theorem \ref{genus g}
and proving semistability of the vertical part by another degeneration argument,
using a comb-like curve with a rational shaft and elliptic teeth contained in fibres, with
the key ingredient being the semistability of the restricted tangent bundle.on the teeth. .\par
\subsection*{Notations and conventions} 
We work over $\C$.
The slope of a vector bundle $E$, i.e. $\deg(E)/\rk(E)$, is denoted $\mu(E)$.
The remainder of $a$ divided by $b$ is denoted $\remainder{a}{b}$. We denote $\remainder {\deg(E)}{\rk(E)}$ 
by $\rho(E)$. The fibre of $E$ at a point $p$ is denoted $E|_p$. It is a $\rk(E)$-dimensional vector space over $\C$.
\subsection*{Acknowledgment} I thank the referee for many incisive comments and corrections.

\section{Some preliminaries}
\subsection{Sheaves and bundles on nodal curves}\label{bundles}
Let $C$ be a connected nodal curve and $E$ a locally free sheaf on $C$.
The fibre of $E$ at $p\in C$, i.e. $E/\m_pE$, is denoted $E|_p$.
For a torsion-free subsheaf $F\subset E$, its double dual $F^{**}$ is reflexive, hence locally free
and we have $F\subset F^{**}\subset E$. Moreover if $F\subset E$ is locally free, then
so is its saturation, i.e. the kernel of the map $E\to (E/F)/(\mathrm{torsion})$.\par
If $F\subset E$ is locally free and saturated, $E/F$ is torsion-free hence locally free
at smooth points, but not necessarily at a node:
\begin{example}\label{node-example}
Let $R=\k[x, y]/(xy)$ be the coordinate ring of a node, with maximal ideal $\m=(x, y)$
and $C=C_x\cup C_y$ the corresponding curve (union of axes).
Then there is an exact sequence
\eqspl{m-resolution}{
0\to \O_C\stackrel{(x, y)}{\to} \O_C\oplus \O_C\stackrel{(-y, x)}{\to}\m\O_C\to 0
}
where $\m\O_C$ is torsion-free but not free.
Thus, unlike in the smooth case, 
it is possible on a nodal curve to have a saturated (cotorsion-free) inclusion
of bundles whose cokernel is not locally free.
\end{example}
\begin{remarks}:\par
\begin{enumerate}\item
 A useful way of looking at a torsion-free sheaf $F$ on $C$, locally near a node $p$,  is
  via the blowup (=normalization) at $p$, denoted $b: \tilde C\to C$, 
 using the device of clutching subspace.
Let $b\inv(p)=\{p_1, p_2\}$ and $\tilde F=b^*F/(\mathrm{torsion})$, which is locally free. 
Then we have
on $C$ an injection of finite colength
$F\to b_*(\tilde F)$.
{\color{blue} This induces a fibre map, which one can check is also injective
\[F|_p\to b_*(\tilde F)|_p=\tilde F|_{p_1}\oplus \tilde F|_{p_2}.\]
The corank of the latter map, i.e. $2\rk(F)-\dim(F|_p)$,  coincides with the rank of the largest local free quotient of $F$.\item
For a saturated injection $j:F\to E$ of locally free sheaves, $ \coker(j)$ is torsion-free but not necessarily
locally free. $\coker(j)$ is locally free iff the induced map $\tilde F\to\tilde E$ is saturated.
}{\color{red}

/****
***/}
\end{enumerate}
\end{remarks}
\subsection{c-semistabilty, openness}
A bundle $A$ on a  connected nodal curve $C$ is said to be (semi)-stable if for every
nontrivial proper torsion-free subsheaf $B\subset A$ one has $\mu(B)\leq\mu(A)$ (resp. $<$
for stable). By the above discussion this is equivalent to assuming the same for
all proper, nontrivial, locally free and saturated, but not necessarily locally cofree,
 subsheaves $B$.\par
In case $C$ is a nodal curve that is the
union of smooth components $C_i$,   a subbundle $B\subset A$ amounts
to a collection of subbundles $B_i\subset A|_{C_i}, \forall i$, together with 
gluing maps
$B_i|_{C_i\cap C_j}\to B_j|_{C_i\cap C_j}, \forall i,j$, compatible with the
gluing maps for $A$, 
and then $\mu(B)=\sum\mu(B_i)$. Note that $B\subset A$ being saturated does not imply
that $B_i\subset A|_{C_i}$ is saturated as the above example shows.
In fact if $B\subset A$ is a saturated subbundle that is not cofree at a node $p$, then
on each component $C_i$ through $p$ the inclusion $B|_{C_i}\to A|_{C_i}$ is unstaturated 
at $p$.
\par
It is immediate from the definition in this case that if $A|_{C_i}$ is semistable
for each $i$ then $A$ is semistable and if moreover $A|_{C_i}$ is stable for
some $i$ then $A$ is stable.\par
$A$ is 
 \emph{cohomologically semistable}
or \emph{c-semistable}  if for every $i\leq\rk(A)$
and every line subbundle (or equivalently, rank-1 subsheaf) $B\subset\wedge^iA$, one has
\[\deg(B)\leq\mu(\wedge^iA)=i\mu(A)=i\deg(A)/\rk(A).\]
By taking determinants, cohomological semistability implies semistability.\par
Also, the above inequality is obviously strict if $i\mu(A)$ is not an integer. $A$ is said to be 
\emph{c-polystable} if is it c-semistable and for each $i$ such that $i\mu(A)$ is an integer,
$\wedge^iA$ is a direct sum of line bundles of the same degree.\par
For a bundle $A$ on a smooth curve we denote by $A_{\max},  A_{\min}$
its first resp. last Harder-Narasimhan bundles (see \S \ref{hn-sec}), i.e. $A_{\max}$ 
(resp. $A_{\min}$) has the largest (resp. smallest) slope,
and largest rank for the given slope, among subbundles (resp. quotient bundles)
of $A$.
The slopes of these are denoted $\mu_{\max}(A), \mu_{\min}(A)$.\par

To some extent this extends to the singular case.
Let $E$ be a vector bundle on a connected nodal curve $C$.
We denote by $\mu_{\max}(E)$ the largest slope of any subsheaf of $E$.
Thus $E$ is semistable iff $\mu_{\max}(E)=\mu(E)$.
If $F\subset E$ has $\mu(F)=\mu_{\max}(E)$ then because $F\subset F^{**}$,
it follows that $F$ is reflexive, hence locally free, and saturated, so that
$E/F$ is torsion-free, hence locally free if $C$ is smooth. Let $F_1\subset F$ be a minimal-rank
subsheaf of the same slope. Thus, every nontrivial proper subsheaf $G\subset F_1$
has slope $\mu(G)<\mu_{\max}(E)$, so $F_1$ is stable. In case $E$ is not semistable, i.e. 
$\mu(F_1)>\mu(E)$, we call $F_1$ a \emphred{stable destabilizer} of $E$. Note that $E/F_1$
cannot be assume locally free, just torsion-free.\par
For a nodal curve $C$, a \emphred{modification} of $C$ is a nodal curve $C'$
with a map $f:C'\to C$ whose non-point fibres are rational trees mapping to nodes of $C$.

The following may be well known but we don't have a reference:
\begin{lem}\label{reduction-lem}
Let $X/B$ be a smooth surface fibred over a smooth curve with smooth general fibre $X_\eta$
and nodal special fibre $X_0$. Let $E$ be a vector bundle on $X$ and $F_\eta$ a subbundle
of $E|_{X_\eta}$. Then there is a base-change $B'\to B$ and a blowup  $X^+\to X_{B'}$
centered on the special fibre and
a subbundle $F_1^+\subset E|_{X^+}$ such that  $F_1^+$ is obtained from $F|_{X^+}$
by a sequence of up modifications centered on components of $X_0^+$ and such that the 
map $F_1^+|_{X^+_0}\to E|_{X^+_0}$ is locally cofree.
\end{lem}
\begin{proof} We may assume the special fibre is a union of smooth components.\par
Step 1: We construct $F_0^+\to E|_{X^+}$ pointwise injective except on nodes of the special fibre.\par
To begin with, as $X$ is smooth, $F_\eta$ extends to a rank-$s$ locally free
subsheaf  $F\subset  E$ on $X$ and by twisting up by components of
the special fibre we may assume the inclusion map $F\to E$
does not vanish on any such component.
If the map$F\to E$ is not injective on some (smooth) component $C$ of the special fibre,
there is a subbundle $K\subset F|_C$ such that $F\subset E$ factors through $F_1$, 
the up modification
of $F$ with respect to $K$. Now we repeat the process for $F_1$.
Ultimately the process ends by the ACC and we have $F_0^+$.\par
Step 2: Now we seek to enlarge $F_0^+$ to obtain a locally free cokernel,
using modifications over nodes only. The process is local
at a node $p$ where the inclusion $F_0^+\to E$ is non-cofree. 
Let $C_1, C_2$ be the branches of the special fibres at $p$.
After a base-change
llocally of the form $t=s^2$ and blowing up $p$, we insert a $D\simeq \P^1$ meeting $C_i$ at $p_i,
i=1,2$. Then on $D$, the map $F^+\to E$ is generically non-injective. Then there is an up modification
$F^+_1$ of $F_0^+$ to which the inclusion $F^+\to E$ extends. 
The process of enlarging $F_0^+$ to $F^+_1$ can be iterated until
the resutling map $F^+_1\to E$ is generically injective on $D$.
Now  consider the maps
\[\wedge^sF_0^+\stackrel{i_1}\to\wedge^sF^+_1\stackrel{i_2}{\to}\wedge^sE, s=\rk(F^+).\]
At $p_1$ and $p_2$, $i_1$ is a non-unit. This means that the vanishing order
of $i_2$ at $p_1, p_2$ on $C_1, C_2$ resp.  is strictly smaller than that of $i_2\circ i_1$. Continuing in this way, or by induction,
we construct a cofree modification of $F^+$.

{\color{red}
/*************
%
%
**********************/
}
\end{proof}
Note it is not claimed nor true in general  that the modification $F^+$ is centered on exceptional components.\par
Our plan for the main results  is to prove (semi)-stability via degeneration of the curve, so
the following result,  saying that (semi)-stability is an open property will be critical:
\begin{lem}\label{openness-lem}[Openness of stability]
Let $\pi:X\to T$ be a family of nodal curves with  with fibres $X_t$  and smooth general fibre. Let $A$ be
a vector bundle on $X$ whose
 restriction on a special fibre  $A_0=A|_{X_0}$ is semi-stable (c-(semi)stable). 
 Then there is a neighborhood $T_0$ of  0
 in $T$ such that for $t\in T_0$, $A_t=A|_{X_t}$ is semi-stable (c-(semi)stable).
\end{lem}
\begin{proof}
We do the semistable case. If the assertion is false then there is a  subbundle 
$B_t\subset A_t$ of degree $>\mu:=i\mu(A_t)$ (in fact, the subbundle of maximal slope
is unique and equals the first Harder-Narasimhan subsheaf). 
After base-changing and blowing up, we obtain a smooth fibred surface $X^+/T^+$
with a map $b:X^+\to X$ such that the $B_t$ glue together to a torsion-free 
subsheaf $B\subset b^*A$ which as in Lemma \ref{reduction-lem} we may assume
injects generically on every component of the special fibre.
This yields a torsion-free sheaf $b_*B\subset A$ hence 
a torsion-free subsheaf $B_0\subset A|_{X_0}$ and, replacing $B_0$ by its double dual
we get a subbundle of $A|_{X_0}$ of slope $>\mu$, contradiction.

{\color{red}
************
*******}
\end{proof}

\subsection{Balanced bundles}
\subsubsection{General case}
A bundle $A$ on a smooth curve  is said to be \emph{balanced} if for every $t$ and a 
general (reduced) effective divisor of degree $t$
one has either $H^0(A(-D_t))=0$ or $H^1(A(-D_t))=0$. 
As mentioned above, balancedness implies
semistability for bundles of integral slope:
	\begin{lem}\label{int-slope-lem}
 	If $A$ is balanced then we have for every subbundle $B\subset A$,
 	\[\mu(B)\leq\ceil{\mu(A)}\]
 	In particular, if $A$ is balanced and  $\mu(A)$ is an integer then $A$ is semistable. 
 	\end{lem}
 	\begin{proof}
 	Adding $1-g$ to both sides, it suffices to prove
 		$\chi(B)/\rk(B)\leq\ceil{\chi(A)/\rk(A)}$.
 	If $t\geq\chi(A)/\rk(A)$ is an integer then $\chi(A(-D_t))\leq 0$, hence $H^0(A(-D_t))=0$.
Therefore $H^0(B(-D_t))=0$, hence $\chi(B(-D_t))\leq 0$, i.e.
 	$t\geq\chi(B)/\rk(B)$.
 	\end{proof}
 	Note the Lemma does not prove that a balanced bundle with integer slope is \emph{cohomologically}
 	semistability.\par
 	\subsubsection{Case of $\P^1$}\label{balance-p1-sec} On $\P^1$, 
 	a bundle $E$ of rank $r$ has the form 
 	\[E=\bigoplus_1^r\O(a_i), a_1\geq...\geq a_r\]
 	and we denote $a_1=\mu_{\max}(E), a_r=\mu_{\min}(E)$.
 	A balanced bundle has the form
 	\[E=a\O(d+1)\oplus(\rk(E)-a)\O(d), \]
 	where $a>0$ and $\deg(E)=\rk(E)d+a$.
 	We call $a$ the 'upper remainder', denoted 
 	$r^+=r^+(\deg(E), \rk(E))$. If $r^+=\rk(E)$, $E$ is said to be perfectly balanced or
 	perfect. Note that this is equivalent to $E$ being semi-stable.  $d+1=\mu_{\max}(E)$ is 
 	also called the upper slope.
 	\subsubsection{Balanced kernels and cokernels}
 	The following elementary results will prove useful in the proof of Theorem \ref{lower} below.
 	\begin{lem}\label{kernel-lem}
 Let
 	\[0\to K\stackrel{\psi}{\to}A\stackrel{\phi}{\to} B\to 0\] 
 	be an exact sequence of vector bundles on $\P^1$ with $A$ balanced.\par
 	(i) If $\phi$ is general for given $A, B$ then $K$ is balanced.\par
 	(ii) If $\psi$ is general for given $K, A$ then $B$ is balanced.
 	\end{lem}
 	\begin{proof} It suffices to prove (i). 
 	We may assume $A, B$ have respective ranks $r>s$.
 	Normalizing, we may assume $A=r_1\O(1)\oplus r_0\O$. Assume $r_0\leq s$ as the contrary case is
similar and 	simpler. Since a surjection $A\to B$ exists, we may assume all line bundle summands  of $B$
 	have degree $\geq 1$, and then we can choose a subbundle $A_0=(s-r_0)\O(1)\oplus r_0\O\subset A$ 
 	with complement $A_1=(r-s)\O(1)$ and a general map $\phi_0: A_0\to B$ such
 	that we have exact
 	\[\exseq{A_0}{B}{\tau}\]
 	where $\tau=\bigoplus k(p_i)$ with $p_i$ distinct points. Clearly, the kernel $K$ of a general map
 	$\bar\phi:A_1\to \tau$ is balanced. 
 	 Because $H^1(\hom(A_1, A_0))=0$,  the map
 	$\Hom(A_1,B)\to \Hom(A_1,\tau)$ is surjective, 
 	hence $\bar\phi$ comes from
 	a map $\phi_1:A_1\to B$ such that $\phi=(\phi_0, \phi_1)$ has kernel $K$.
 	\end{proof}
 	Sometimes the hypothesis that the middle bundle be balanced can be relaxed:
 	\begin{lem}\label{coker-lem}
 	(i) Let $A, B$ be bundles of respective ranks $r+1, r$  on $\P^1$ such that
 	$A^*\otimes B$ is globally generated.  Then there is a surjection $A\to B$
 	with  kernel $\O(k)$ where $ k:=\deg(A)-\deg(B)$.\par
 	(ii) Let $A, k$ be such that 
 	\[\mu_{\max}(A)\leq (\deg(A)-k)/(\rk(A)-1).\] 
 	Then the cokernel
 	of a general map $\O(k)\to A$ is balanced.
 	\end{lem}
 	\begin{proof}
 	(i) Trivial.\par
 	(ii) Let $B$ be the unique balanced bundle of rank $\rk(A)-1$ and degree $\deg(A)-k$
 	hence lower degree $[(\deg(A)-k)/(\rk(A)-1)]$.
 	Then our numerical hypothesis implies that $A^*\otimes B$ is globally generated, hence
 	there exists a surjection $A\to B$ with kernel $\O(k)$. Because balancedness is an open
 	property, it follows that the cokernel of a general map $\O(k)\to A$ is balanced
 	(and isomorphic to $B$).
 	\end{proof}
 	Note that the condition $A^*\otimes B$ globally generated means
 	$\mu_{\max}(A)\leq \mu_{\min}(B)$, i.e. (max degree of line subbundle of $A$)$\leq$
 	(min degree of line bundle quotient of $B$).\par
 	For application (see Theorem \ref{lower}), we require an extension of Lemma \ref{kernel-lem} 
 	to some embedded situations where the map 
 	$\O(k)\to A$ is
 	not general and $A$ is unbalanced. To this end, let 
 	\[A=\bigoplus\limits_1^{r+1}\O(a_i), a_1\geq...\geq a_{r+1}, r>0\]
 	be a completely decomposable bundle on $\P^m$, and 
 	suppose we have exact
 	\eqspl{b-on-p}{\exseq{\O(k)}{A}{B}
 	}
 	where the map $\O(k)\to A$ is general, given by general homogeneous polynomials
 	on $\P^m$. Note that $B$ is locally free outside of a codimension $r+1$ subset.
Now	let $C\subset\P^m$
 	be a general rational curve of degree $e$. Then  $B|_C$ is locally free and we get exact on $C$
 	\eqspl{b-on-c}{\exseq{\O(ek)}{\bigoplus\O(ea_i)}{B|_C}.
 	}
 	If we assume $C$ is projectively normal, or equivalently that $e\leq m$, we have that the left map in
 	\eqref{b-on-c} is a general map, given by general
 	homogeneous polynomials on $C$, and it it follows from Lemma \ref{coker-lem} that
 	$B|_C$ is balanced hence perfect if $r|e$. We now extend this observation  to 
 	curves of higher degrees $e$.
 	\begin{lem}\label{g-on-c-lem}
 	Notations as above assume $C\subset\P^m$ is a general rational curve of degree $e$
 	and that\par
 	(i)  $a_1\leq (\sum a_i-k)/r$; \par
 	(ii) $1\leq r\leq m$. \par
 	Then $B|_C$ is balanced.
 	\end{lem}
 	Note that when $r<m$, $B$ is not locally free but $B|_C$ is.
 	\begin{proof}
 	Write $e=ur+v, v\in [0,r)$ and consider a general rational chain
 	\[C'=C_0\cup_{p_1}C_1\cup_{p_2}...\cup_{p_u}C_u\]
 	where $C_0$ has degree $v$ and $C_i$ has degree $r$ for $i>0$.
 	Then $B|_{C_0}$ is balanced while for $i>0$, $B|_{C_i}$ is perfect.
 	Hence $B|_{C'}$ is balanced, as is the restriction on a general smoothing of $C'$.
 	\end{proof}
 	\subsubsection{Upper and lower subbundles}
 Notations as above.
 	The subbundle $a\O(d+1)$  is uniquely determined and called the \emphred{upper subbundle}
 	$E^+$.
 	A complementary to the upper is call a \emphred{lower} subbundle. It is not unique
 	and in fact it can be chosen so  its fibre at two  given points $p ,q$ is any given pair of 
 	complementary subspaces $L_p, L_q$ to 
 	the fibres $E^+|_p, E^+|_q$ of the upper
 	subbundle. This assertion  follows easily from the fact that the kernel of the surjection
 	\[E\to E|_p/L_p\oplus E|_q/L_q,\]
 	i.e. the down modification corresponding to $L_p, L_q$-see the next subsection-  is balanced 
 	and has the form $(\rk(E)-a)\O\oplus a\O(-1)$,
 	hence has $h^1=0$.
 	\subsubsection{Perfect curves}
 	A rational curve $C\subset X$ is said to be perfect if its normal bundle is.
 	If $n$ is understood, an integer $e$ is perfect if $2(e-1)\equiv 0\mod n-1$.
 	A general rational curve of degree $e$ in $\P^n$ is perfect if $e$ is.
 	\subsection{Bundle modifications}
 	A \emphred{down modification} of a bundle $E$ on a variety $X$ is the (locally free)
 	kernel of a surjection $E\to\tau$ where $\tau$ is a locally free $\O_D$-module
 	for some (possibly nonreduced) Cartier divisor $D$ on $X$. An \emphred{up modification} of $E$ 
 	is an injection $E\to G$ of bundles whose cokernel is a locally free $\O_D$-module; or
 	equivalently, the dual of a down modification of the dual $\check E$.\par
 	For a bundle $E$ on a curve, a \emphred{general down modification} of $E$ is the kernel of a torsion quotient
 	$E\to \tau$  such that $\tau$ is supported on a general collection of points and the map
 	$E\to \tau$ is general.
 	Specifying such a modification is equivalent to specifying some 
 	distinct general points $p_1,...,p_t$ plus general 
 	subspaces $V_i\subset E|_{p_i}, i=1,...,t$ and then $\tau=\sum (E|_{p_i}/V_i)k(p_i)$.
 	Note that general modifications make sense even if the curve is reducible, in which case it
 	is assumed that each point $p_i$ us general in some component.
 	It is easy to see that a general modification of a balanced bundle on $\P^1$ is balanced.
\par
A \emphred{general up modification} of $E$ is by definition an injection of bundles  $j:E\to G$
whose cokernel is supported on a general $t$-tuple $p_1,...,p_t$ and such that the kernel of
$E|_{p_i}\to G|_{p_i}$ is a general subspace, for each $i$; equivalently, a general up modification
of $E$ is the dual of a general down modification of the dual $E^{*}$.\par
A modification as above is \emph{simple} if its cokernel has length $\leq 1$ locally at each point.
 \par 
A bundle is said to be \emph{hyper-stable} if its general down  or up modification is stable.
Ditto for semistable and for the c-versions.\par
It is easy to see that a general (up or down) simple modification of a balanced bundle is balanced.
However a general modification of a stable bundle need not be semistable.
\subsection{Quot schemes}\label{quot-sec}
See \cite{grothendieck-hilbert},  \cite{sernesi}. Given a bundle $E$ on a scheme $X$, 
quotients, or equivalently subbundles $F\subset E$ are parametrized by a 
Grassmannian or Quot  scheme $\mathrm{Quot}_E$ 
whose respective  tangent and obstruction spaces are $H^0(\hom(F, E/F))$
and $H^1(\hom(F, E/F))$. If $X$ is a curve of arithmetic genus $g$, 
and $F, E$ have respective ranks $s, r$, we have
\[\chi(\hom(F, E/F))=(h^0-h^1)(\hom(F, E/F))=-r\deg(F)+s\deg(E)+s(r-s)(1-g).\]
 \subsection{HN filtration}	\label{hn-sec}
 A bundle $E$ admits a uniquely-determined increasing filtration 
 	called the  Harder- Narasimhan filtration
 	\[(E_\bullet)=(E_0=(0)\subsetneq E_1\subsetneq E_2\subsetneq...
 	\subsetneq E_k=E)\] such that each $E_i/E_{i-1}$ is semistable and
 	$\mu(E_i/E_{i-1})>\mu(E_{i+1}/E_i)$. In particular, $E_1$ is called the first HN subbundle of $E$
 	and is characterized by having maximal slope, 
 	denoted $\mu_{\max}(E)$,
 	and maximal rank for its slope, among subsheaves of $E$.
 	A \emphred{stable destabilizer} is a minimal-rank 
 	nonzero subsheaf $G\subset E_1$ of slope $\mu(G)=\mu(E_1)$. 
 	
 	/****************
 \subsection{Fangs and their hypersurfaces}\label{fangs}
 See \cite{caudatenormal} for details.\par
 	Given a linear $A_1=\P^{n-m-1}\subset\P^n$,
 	let $\pi_1:Z_1=B_{A_1}\P^n\to\P^m$ be the resolved projection centered
 in $A_1$
 	with blowdown map $b_1:Z_1\to\P^n$ and exceptional divisor $E_1=A_1\times\P^m$. 
 	Let $\pi_2:Z_2=B_{A_2}\P^n\to\P^{n-m-1}$ be the analogous
 	object, based on blowing up $A_2=\P^m$, with blowdown map $b_2$
 	and exceptional divisor $E_2$. Setting $E=\P^m\times\P^{n-m-1}$,  
 	an arbitrary pair of isomorphisms
 	\[\phi_1:\P^m\to A_2, \phi_2:A_1\to\P^{n-m-1}\] yield 
 	isomorphisms \[E\to E_1, E\to E_2,\] whence 
a normal-crossing
 	variety
 	\[Z_0=Z_1\cup_EZ_2\]
 	called a \emph{generalized fan} or \emph{fang} of type $(n,m)$.\par
 	A flat morphism $\cZ\to B$ is called a \emph{relative fang} of
 	type $(n,m)$ if each fibre is either $\P^n$ or a fang of type $(n,m)$.
 	A standard way to construct a relative fang is to blow up
 	the subvariety $\P^{n-m-1}\times 0$ in $\P^n\times\A^1$.
 	There $Z_1$ and $Z_2$ are, respectively, the birational transform of $\P^n\times 0$
 	and the exceptional divisor, and both have projective bundle structure:
 	\[Z_1=\P(G_1), Z_2=\P(G_2)\]
 	where
 	\[G_1=\O_{\P^m}(1)\oplus (n-m)\O_{\P^m}, G_2=\O(1)_{\P^{n-m-1}}\oplus (m+1)\O_{\P^{n-m-1}}.\]
 	We will denote the latter as $Z_1=\P_{\P^m}(1,0^{n-m}), Z_2=\P_{\P^n-m-1}(1, 0^{m+1})$.
 Also each $Z_i$ is endowed with the $\O(1)$ induced by from $G_i$ which is
 	also the pullback of $\O_{\P^n}(1)$ by the blowdown map $b_i:Z_i\to\P^n$.
 	We denote the exceptional divisor of $b_i$ by $E_i$, and we also have
 	\[E_1=\P((n-m)\O_{\P^m}), E_2=\P((m+1)\O_{\P^{n-m-1}}),\]
 	hence of course
 	\[E=E_1\simeq E_2\simeq \P^m\times\P^{n-m-1}.\]\par
 The case $m=0$ or $m-n-1$ reduces to the fan construction,
 	Now assume $0<m<n-1$ and $e<d$. Then
 	the linear system $|dH-eZ_2|$ on $\cZ$, 
 	where $H$ is the pullback of a hyperplane in
 	$\P^n$, restricts as follows.
 	\begin{itemize}\item on the general fibre, to $|dH|$;\item
 	on $Z_1$,  to $|db_1^*H_{\P^n}-eE_1|$, i.e. the birational transform on $Z_1$
 	of the system of hypersurfaces of degree $d$ on $\P^n$ 
 	with multiplicity $e$ on $\P^{n-m-1}$;\item
 	on $Z_2$ to $|db_2^*H_{\P^n}-(d-e)E_2|$;\item on 
 	$E=\P^m\times\P^{n-m-1}$ to the linear system of hypersurfaces
 	bidegree $(e, d-e)$.
 	\end{itemize}
 \begin{example}\label{curve-degree} Taking $d=1, e=0$  shows that a curve of degree $k$ in 
 $\P^n$ can specialize to a curve $C_1\cup C_2\subset Z_0$ where $C_1\subset Z_1$ maps to
 a curve $b_1(C_1)$ of degree $k_1$ in $\P^n$ while $C_2$ 
 a curve
 $b_2(C_2)\subset\P^n$ of degree $k-k_1+\l$ meeting te blown-up $\P^m$ in $\l$ points (hence $C_2$
 maps to a curve $\pi_2(C_2)$ of degree
 $k-k_1$ in the base $\P^{n-m-1}$). 
  Here $\l=C_1.E=C_2.E$ and $C_1\cap E=C_2\cap E$ consists of $\l$ points.
  \par Similarly for the dual case $d=e=1$.\par
  In the case $m=0, e=d-1$ we obtain the quasi-cone degeneration used previously.
 	\end{example}
 	The foregoing construction may obviously be extended to the case of
 	more than 2 components but we don't need this.
 \section{More preliminaries: normal bundles of reducible curves}
 Let $X$ be a smooth quasi-projective variety of dimension $n$ and let $C_1, C_2$ be
 locally complete intersection (lci)
 curves on $X$ meeting  at a unique point $p$ smooth on both and with distinct tangent directions. 
 Then $C=C_1\cup C_2$ is a lci curve
 hence has a locally free normal bundle $N$.\begin{remarks}
 Some basic properties of   $N$ are:
\begin{enumerate}\item  If $C$
 is contained in a smooth surface germ $S$ then we have an exact sequence
 \[\exseq{\O_S(C)\otimes\O_C}{N}{N_{S/X}\otimes\O_C}.\]
\item Let $\tau=T^1(C)=\Ext^1(\Omega_C, \O_C)$, a length-1 skyscraper sheaf at $p$.
 Then we have exact
 \[\exseq{N^0}{N}{\tau}\]
 where $N^0$ is the (non-locally-free)  'locally trivial normal sheaf' , image of
 the natural map $T_X\to N$, which corresponds
 to deformations of the triple $(C_1, C_2, p)$ and which fits in the following exact sequence
 in which $N_i=N_{C_i/X}, N_{S,p}=T_{p,X}/(T_{p,C_1}+T_{p,c_2})$;
 \[\exseq{N^0}{N_{1}\oplus N_{2}}{N_{S,p}}.\]
 In particular, $N^0|_{C_i}=N_{C_i/X}$. Thus $N|_{C_1}$ is an elementary up
 modification of $N_{C_1/X}$, i.e. there is exact
 \[\exseq{N_{C_1/X}}{N|_{C_1}}{\tau}\]
 inducing
 \[0\to T_{p, C_2}\to N_{C_1/X}|_p\to N|_p\to\tau\to 0.\]
 Hence we have exact
 \[\exseq{N_{S}|_p}{N|_p}{\tau}.\]\par
\item Now let $v_{1,1}\in N_1|_p$ generate $T_{C_2, p}$, and extend and enlarge $v_1$ to a 
 local basis $(v_.)$ of
 $N_1$. Let $w_{1,1}\in N$ be a local  lift of $\tau$
 Then $w_{1,1}$ can be extended to a  local basis of $N$ so that  the inclusion map $N_1\to N|_{C_1}$
 has matrix $\diag(t, 1,...,1)$ therefore a local section of $N_1$ of the form
 $v_{1,1}+t(v_{2, 1}+...)$ maps to a section of $N|_{C_1}$ to the form $t(w_{1,1}+w_{2,1}+...)$
 and the latter saturates to the line subbundle  generated by $w_{1,1}+w_{2,1}+...$.
 Ditto for $v_{2,i}, w_{2,i}$ etc. Now note that $w_{1,i} $ and $w_{2,i}$ correspond to 
 $\O_S(C)|_C$ where $S$ as above is a germ of a smooth surface containing $C$.
 Therefore we may assume $w_{1,1}$ and $w_{2,1}$ agree at $p$ so they generate
 the same local line  subbundle of $N$ that we call a \emph{singular} subbundle.
 It is in general not unique though in some cases it is unique by virtue of being the
 maximal slope subbundle.
 \par
 \item  Note that $t_{S,p}:=N_{C_1/S}|_p$ coincides with the image of $T_{C_2,p}$ in $N_1|_p$
 and maps isomorphically onto $\tau_p$. Therefore we have canonically
 \[N_1|_p=t_{S,p}\oplus N_{S,p}\]
 and $N_{S,p}\subset N|_p$ coincides with the image of $N_1|_p$ and that of
 $N_2|_p$.\par
 Similarly, we have
 \[N|_p=\tau\oplus N_{S,p}\]
 and we call the summands to \emph{singular} resp. \emph{regular} subspaces.
 Thus the following sequence is exact
 \[0\to t_{S, p}\to N_1|_p\to N|_p\to\tau\to 0.\]
 Yet another way to think of this setup is via local equations: equations for $C$
 have the form $x_1x_2, x_3,...,x_n$ where $x_3,...,x_n$ are equations for $S$
 and $x_1, x_2$ are respective equations for $C_2, C_1$ on $S$. Then 
 $x_3,...,x_n$ correspond dually to $N_{S,p}$ while $x_2$ corresponds
 to $t_{S,p}$ on $N_1|_p$.\item
 It follows from the above discussion that the gluing of $N|_{C_1}$ and $N|_{C_2}$
 at $p$ is not general. This fact is partly responsible for the difficulty in proving
 stability of normal bundles via degeneration.
 \end{enumerate}
 \end{remarks}
 \begin{examples}\label{examples}
  We use notations as above and also, for $\P^1$,  set $\O(a_\bullet)=\bigoplus\O(a_\bullet)$,
 and for $C_1\cup_p C_2$, set $\O(a\cup b)=\O_{C_1}(a)\cup\O_{C_2}(b)$.
 Also, a bundle $E$ on $\P^1$ is \emph{balanced} if $E\simeq \O(*(a+1), *a)$ for some $a$.
 \begin{enumerate}\item
 A conic $C_1$ in $\P^3$ has normal bundle $\O(4,2)$ while a line $C_2$
 has normal bundle $2\O(1)$. If $C=C_1\cup_pC_2$ Then
 \[N|_{C_1}=\O(4,3), N|_{C_2}=\O(2,1)\]
 with $w_{1,1}$ generating an $\O_{C_1}(3)$ and $w_{2,1}$ generating the $\O_{C_2}(2)$.
 Thus these glue together to a line subbundle $\O(3,2)$ while the other basis elements generate 
 a $\O(4,1)$. Thus
 \[N=\O(4\cup 1)\oplus\O(3\cup 2)\]
 which smooths out out to $2\O(5)$ on a twisted cubic. Note that the $\O_{C_1}(3)$ and $\O_{C_2}(1)$
 are not unique, however $\O(4\cup 1)$ and $\O(3\cup 2)$ are unique because
 there is a unique $\O_{C_1}(3)$ that glues at $p$ to the unique $\O_{C_2}(2)$ and ditto for
 $\O_{C_2}(1)$.  \par
 \item Similarly, if $C_1\subset\P^n\subset\P^{n+1}$ is a rational normal curve and $C_2\subset\P^{n+1}$
 is a general line through $p\in C_1$ then
 \[N_1=(n-1)\O(n+2)\oplus\O(n+1), N_2=\O(2)\oplus (n-1)\O(1)\]
 then the $\O(n+1)$ glues to the $\O(2)$ so that
 \[N=\O((n+1)\cup 2)\oplus\bigoplus(n-1)\O((n+2)\cup  1)\]
 which smooths out to $n\O(n+3)$.\par
 \item If $C_1, C_2\subset\P^3$ are conics, we have
 \[N|_{C_i}=\O(4,3).\]
 There are canonical choices for $\O_{C_1}(3)|_p=\O_{C_2}(3)|_p$
 while $\O_{C_1)(4)|_p\neq \O(c_2}(4)|_p$. However the family of all $\O_{C_1}(3)$
 subsheave fills up an open subset of $N|_p$ containing $\O_{C_2}(4)|_p$
 and vice versa. Therefore we have
 \[N=\O(4\cup 3)\oplus\O(3\cup 4)\]
 which smooths out to $2\O(7)$ on a rational curve of degree 4.
 \item Similarly if $C_1\subset\P^3$ is a rational curve of degree $e$ with normal bundle 
 $2\O(2e-1)$ and $C_2$ is a conic, we have
 \[N_C=\O(2e\cup 3)\oplus\O(2e-1\cup 4)\]
 which smooths out to $2\O(2e+3)$ on a rational curve of degree $e+2$.
 \item Inductively, this together with the initial cases $e=3, 4$ 
 recovers the well-known fact that a general rational curve of degree
 $e$ in $\P^3$ has normal bundle $2\O(2e-1)$.\item
 Let $n=4$. For  $e=5$ we use $C_1, C_2$ of degrees 3,2  which leads to
 \[N|_{C_1}=\O(5,5,4), N_{C_2}=\O(4,3,2).\]
 There are 1-parameter families of $\O_{C_1}(5)$ and $\O_{C_2}(3)$ so we can glue 1 pair
 to get a $\O(5\cup 3)$ subbundle. Then since there are 4 and 5-dimensional families of $\O_{C_1}(4)$
 resp $\O_{C_2}(2)$ we can get complementary subbundles $\O(5\cup 2), \O(4\cup 4)$.
 Altogether this smooths out to $2\O(8)\oplus\O(7)$ on $e=5$.\par
 Next for $e=6,7$ a simpler argument with $C_2$ a line, so that $N|_{C_2}=\O(2)+2\O(1)$
 yields $N=\O(8\cup 2)\oplus\O(8\cup 1, N=\O(9\cup 2\oplus 2\O(10\cup 1)$
 which smooth out to $\O(10+2\O(9)$ resp, $3\O(11)$.\par
 For general $e=3a+b>7$ we use induction and attachment of conic rep. line resp. line 
 as above depending on 
 $b=2, 0, 1 \mod 3$ leading eventually   to $2\O(5a+3)\oplus\O(5a+2), \O(5a)+2\O(5a-1), 
 3\O(5a+1$ resp. on a general degree-$e$ rational curve.
 \item Generally, suppose $N_1$ is balanced and has coremainder $k>0$ and $C_2$ is a rational normal curve in
 a $\P^k\subset\P^n$. Then
 \[N=(k-1)\O(a\cup k+2)\oplus\O(a+1\cup k+1)\oplus (n-k-1) \O(a+1\cup k).\]
 This smooths out to a balanced bundle of remainder $k$.
 \par
 \item Assume $C_1$ has degree $e_1\equiv 1\mod n-1$ so $N_1$ is perfect, i.e. 
 \[N_1=(n-1)\O(a), a=e_1+(2e_1-2)/(n-1),\]
 while $C_2$  is a rational normal curve in $\P^{n-1}\subset\P^n$,
  Then 
  \[N|_{C_1}=\O(a+1)\oplus(n-2)\O(a), N|_{C_2}=(n-2)\O(n+1)\oplus\O(n)\] where
  $\O_{C_1}(a+1)$ glues to $\O_{C_2}(n)$ so 
  \[N=\O((a+1)\cup n)\oplus (n-2)\O(a\cup (n+1))\]
  which smooths out to a perfect bundle on $C$ of degree $e=e_1+n-1$.
  \item (Line attachment) Given $e$ of remainder $k>0$, we can take $C_2$ as a union of $e_2=n-1-k$ lines
  with $C_1$ general of degree $e_1=e-e_2$ meeting each one. Then $N$ is a direct sum
  of $e_2$ line bundles of the form $\O((a+1)\cup 2\cup 1...)$ and permutations
  and ones of the form $\O((a_1)\cup 1\cup...)$ and again smooths out to a balanced bundle.
 \end{enumerate}
 \end{examples}
 Building on the Examples, we sketch a short  proof of the following well-known result; this proof, and the 
 specialization involved, will be used in the sequel.
 \begin{prop}\label{balanced}
 A general rational curves of degree $e\geq n\geq 3$ in $\P^n$ has balanced normal bundle.
 \end{prop}
 \begin{proof} Case 1: $e=n$ (rational normal curve): discussed above.\par
 {\color{blue}
 Case 2: $n<e\leq 2n-1$. Write $e=n+\l, m=n-\l\geq 1$. Consider a reducible curve
 \[C_0=C_m\cup\bigcup\limits_1^{2\l}M_i\]
 where $C_m$ is a rational normal curve in $\P^m\subset\P^n$ and the $M_i$ are general
 1-secant lines of $C_m$ in $\P^n$. Then we have
 \[N_{C_m/\P^n}=(m-1)\O(m+2)\oplus (n-m)\O(m),\]
 \[N_{C_0/\P^n}|_{C_m}=(n-1)\O(n-\l+2), N_{C_0/\P^n}|_{M_i}=\O(2)\oplus (n-2)\O(1).\]
Hence
 $N_{C_0/\P^n}$ is a direct sum of line bundles of total degree $n+\l+2$ and $n+\l+3$ so 
 smooths out to a balanced bundle.\par
 Case 3:$e\geq 2n$. Then consider a curve of the form
 \[C_0=C_1\cup_p C_{n-1}\]
 where $C_1, C_{n-1}$ are general of degrees $e-n+1, n-1$. Using induction and the
 rational normal case it is nor hard to see that $N_{C_0/\P^n}$ smooths out to a balanced bundle. 
 
 }
 {\color{red}
/*********** 
}.

 {\color{red}
 *********
 *************************/
 }
 \end{proof}
 \begin{rem}\label{normal-degenerate-rem} Note that for $e<n$ the normal bundle $N_{C_e/\P^n}$ has the form
  $(e-1)\O(e+2)\oplus (n-e)\O(e)$
 so it is unbalanced; the splitting 
 is non-canonical. But if we attach $k$ general lines to $C_e$ with $n-e\leq k\leq 2(n-e)$ then the
 normal bundle of the union restricts to a balanced bundle on $C_e$.
 \end{rem}
 	
 \section{Parallel transport and holonomy}
\subsection{General case} Let  $C$ be a smooth rational curve and   $N$ a vector   bundle on $C$. The $N$
 admits a canonical increasing filtration $F.N$ called the Harder-Narasimhan filtration 
 with the property that
 \[\gr_iN=F_iN/F_{i-1}N\simeq V_i(e_i):= V_i\otimes\O(e_i).\]
 For points
  $p\neq q\in C$, we have  restriction isomorphisms, canonically  defined up to up to scalar factor
  \[\gr_iN|_p\stackrel{\sim}{\leftarrow} H^0(\gr_iN(-e_i))\simto \gr_iN|_q.\]
  Now assume that we are given filtered  vector spaces $W_p, W_q$, 
  that remain fixed  as the pair $(C, N)$ moves,  together with
  isomorphisms $\gr_iW_p\simeq \gr_iN|_p, \gr_iW_q\simeq \gr_iN|_q, \forall i$. 
 Then the resulting isomorphism $\gr_\bullet W_p\simto \gr_\bullet W_q$, which depends on $C$,
 is well defined up to a 'graded scalar', i.e. scalar on each summand. The induced map
 \[\tau_{C,N,  p,q}:\prod\P\gr_iW_p\to\prod\P\gr_iW_q\] 
 is called a (graded) \emph{parallel transport map}
 (associated to this data). \par
 Given a cycle (closed $n$-gon) of smooth rational curves
 \[C=C_1\cup_{p_1}C_2\cup...\cup_{p_{n-1}}C_n, C_n\cap C_1=p_n,\]
 assume for simplicity that $N|_{C_i}$ is perfect for each $i$. Then the \emph{holonomy}
 map based at a node $p_1$
 is by definition the composite map 
 \[\eta_{C, N, p_1}=\tau_{p_n, p_1}\circ\tau_{p_{n-1}, p_n}\circ...\tau_{p_1, p_2}\in \P(\gl(N_{p_1})).\]
\subsection{Normal  and restricted bundles}  One case of these constructions of interest to us is the following. Let $X$ be a smooth variety, $p\neq q$
 points on $X$, and $D_p, D_q$ germs of smooth divisors ar $p$ resp. $q$. Set
 \[W_p=T_pD_p, W_q=T_qD_q.\]
 Now let $C\subset X$ be a rational curve, smooth and transverse to $D_p, D_q$
 at $p,q$, with normal bundle $N$ filtered as above. 
 The natural maps $T_pX\to N|_p, T_qX\to N|_q$ yield canonical isomorphisms
 \[W_p\simeq N|_p, W_q\simeq N|_q.\]
 These endow $W_p, W_q$ with canonical HN filtrations.
 \begin{defn}
 The map
 \[\tau_{C, p, q}: \prod\P\gr_iW_p\to \prod\gr_iW_q\]
 defined by the above data is called the associated \emph{(normal) parallel transport}.\par
 Similarly if $C=C_1\cup_{p,q}C_2$ where $C_1, C_2$ are perfect rational curves,
 the map
  \[\eta_{C, p}\in\P(\gl(N_p)),\]
 \[\eta_{C, p}=\tau_{C_2, q, p}\circ\tau_{C_1, p,q}=\tau_{C_2, p,q}\inv\circ\tau_{C_1, p,q}.\]
 is called the associated \emph{holonomy}.
 \end{defn}
 Another case of interest is that of the restriction of a bundle such as the tangent bundle
 from the ambient variety. Thus if $E$ is a bundle on $X$ such that the restriction
 $E|_C$ is perfect, we get a similar parallel transport map, defined up to scalar multiple
 \[\tau_{C, E, p,q}:E|_p\to E|_q.
 \]
 In case $X=\P^n, E=T_X$, the parallel transport map is distinct from a map $T_{X, p}\to T_{X, q}$
 induced by a projective transformation mapping $p$ to $q$.
\subsection{Quasi generality} 
 Our main interest is in the dependence of the parallel transport map on $C$ with fixed $p,q$
 and fixed filtered vector spaces $W_p, W_q$; in particular
 how general it is for general $C$, fixing $p,q$. The latter
 question is independent of choices. One consequence of generality
 is that we can glue two such curves through $p,q$, thus obtaining a genus-1
 curve whose holonomy is a general endomorphism, hence   diagonalizable with distinct eigenvalues.
 \par
 To make this precise first note the following. Suppose given a family $\{C_t:t\in T\}$ of rational curves
 passing through fixed points $p\neq q\in X$ and transverse to divisors $D_p, D_q$,
 and such that the respective normal bundles $N_t$ all have the same type
 (e.g. are balanced or even perfect). Then as above, we get 
 a family of maps, depending on $t$,
 \[\tau_{C_t, p, q}=\tau_{C_t, N_t, p,q}\in \Hom(W_p, W_q).\]
 Up to composing with  automorphisms of $W_p, W_q$ that are fixed independent of $t$, this 
 family is independent of choices. Therefore it makes sense to consider the size of this family,
 in particular whether it is dense. Typically, we are able to prove a slightly weaker property
 called quasi generality, 
 which leads us to the following.

 \begin{defn} (i) Let $V, W$ be $n$-dimensional vector spaces and
 consider an irreducible subvariety $T\subset\P(\Hom(V, W))$ of linear maps with general member  $\tau$. 
 Then $T$ or $\tau$ is said to be a \emph{quasi-general} (qg)
 or \emph{ quasi-mixing} if it is 'projetively generically $n$-transitive', i.e. for some or equivalently, almost any) basis 
 $(v.)=(v_1,...,v_n)$ of $V$,  the $n$-tuple $([\tau v_1], ...,[\tau v_n])\in\P(W)^n$,
 as $\tau$ varies,  is a  general point;
 i.e. the  map $T\to\P(W)^n$ determined by evaluation on  $(v.)$ is dominant. \par 
 A graded map of graded vector spaces is quasi general if each graded piece is.\par
 A filtered map of filtered vector spaces is quasi general if the associated graded is.\par
(ii) Suppose given a family  $(\mathcal C, \mathcal N)=\{(C_t, N_t):t\in T\}$ of pairs 
of the form  (rational curve, balanced bundle), 
together with disjoint sections $p,q$ and rigidifications at those,
 i.e. trivializations
$\mathcal N|_p\simeq W_p\O_T, \mathcal N|_q\simeq W_q\O_T$.
The family, or its general member $(C_t, N_t)$  is said to be \emph{quasi-mixing} 
if the family of associated parallel transport maps $\{\tau_{C_t, N_t, p,q}\}$ is a quasi-general in 
$\P(\Hom(W_p, W_q))$. 
\par
(iii) A rational curve $C$ on a variety $X$ is said to be \emph{quasi-mixing}
if for general points $p,q\in C$. the pair $(C, N_{C/X})$ together with $p,q$
are members of a quasi-mixing family.
\end{defn}
\begin{remarks}\begin{enumerate}\item
For a  quasi-projective subvariety  $T\subset\P(\Hom (V, W))$
with general element $\tau$,
$T$ will be transverse to a general translate 
by $\gl(V)\times\gl(W)$
of the set of diagonal matrices with respect to some bases of $V, W$, therefore
 the property of being quasi-general is equivalent to $\dim(T)\geq n(n-1)$
 (i.e. the affine cone on $T$ should have dimension at least  $n^2-n+1$. Hence
 $T$ is quasi-general iff for \emph{almost every} basis $(v.)$, $([\tau v.])$ is a general $n$-tuple.
 For example the set of  upper triangular matrices in $\PGL_2$ is quasi-general though it
 fixes $[1,0]$.
\item Consequently $\{\tau\}$ is quasi-general iff for some $v_1,...,v_k\in V$ independent and
some $w_1,...,w_{n-k}\in W$ independent, the n-tuple 
\[([\tau(v_1)],...,[\tau(v_k)], [\tau\inv(w_1)],...
[\tau\inv(w_{n-k})])\in \P(W)^k\times\P(V)^{n-k}\] is general.\item
Quasi-general is weaker than general or mixing. A subvariety  $T$ is dense in $\P(\Hom(V, W))$, i.e. its general member $\tau$ is general, 
iff for one (or equivalently, almost every) $n+1$ (rather than $n$) tuple $([v_0],..,[v_n])\in\P(V)^{n+1}$
in general position, the $(n+1)$-tuple $([\tau(v.)])\in\P(W)^{n+1}$ is general, i.e.
the family $\{([\tau(v.)]):\tau \in T\}\subset (\P(W))^{n+1}$ is dense.
Thus if $\tau$  is quasi-general  and $[\tau(v_0)]$ is general
independently of $([\tau(v_1)],... ,[\tau(v_n)])$ then $\tau$ is general. 
Thus e.g. for $n=2$, i.e. $\P^1$, quasi-general (resp. general) means
generically 2 (resp. 3) transitive on  $\P^1$. For general $n$, $T$ is general iff
for \emph{every} $(n+1)$-tuple in general position in $\P(V)^{n+1}$, its $T$- orbit
is dense in $\P(W)^{n+1}$. \item
If $\tau_1\in\Hom(V_1, V_2)$ and $\tau_2\in\Hom(V_2, V_3)$ are quasi-general
then $\tau_2\circ\tau_1$ is general. In particular, if $V_3=V_1$ then
$\tau_2\circ\tau_1\in \gl(V_1)$ is diagonalizable with distinct eigenvalues.
Consequently if $C=C_1\cup_{p,q}C_2$ is a fish curve such that for some bundle $E$,
e.g. the normal bundle $N_C$,
$\tau_{C_1, E, p,q}$ and $\tau_{C_2, E, q, p}$ are quasi general then the holonomy
at $p$, i.e. 
$\eta_{C, E, p}=\tau_{C_2, E, q, p}\circ\tau_{C_1, E,  p,q}$ is diagonalizable with distinct, general eigenvalues. 
Consequently $E$ is a direct sum of distinct, general line bundles as are
its exterior powers $\wedge^E$.
\item If $N$ is a perfect bundle of rank $r$ and $N'\subset N$ is a general simple modification at $\leq r$ 
general points, then 
$N'$ is balanced and admits a  line bundle decomposition compatible
with that of $N$ so   $\tau_{N',p.q}$
can be identified with  $\tau_{N, p,q}$. Hence,  if $(C, N)$ is quasi mixing and $N'\subset N$ is a general simple modification at $\leq r$ points 
then $N'$ is quasi mixing.\item
The quasi mixing property does not depend on the choice of divisors $D_p, D_q$. \item
In applications of quasi-generality, $N$ is balanced and usually perfect.
\end{enumerate}
\end{remarks}

\subsection{Some (quasi-)mixing curves}
We begin by proving quasi- mixing and mixing for a perfect modification of the normal bundle 
of a rational  curve in projective space:
\begin{lem}\label{transport-Pn}
 Let $p\neq q$ be points in $\P^n, n\geq 3$ and
let $C$ be a general rational curve through $p,q$  of degree $e\geq 2n-1$
and let $k\geq 0$ be   such that
 \[k\equiv 2e-2 \mod n-1.\] Let $N'$ be a general simple
down modification of $N_{C/\P^n}$ of colength $k$.  
Then\par
(i) $N'$ is perfect and quasi mixing, i.e. as $(C, N')$ vary
fixing general points  $p,q\in C$  independent of the modification points, the family $\{\tau_{C, N'}\}$ is quasi-
general.\par
(ii) If moreover $e\geq 4n-3$ then the family is mixing, i.e. $\{\tau_{C, N'}\}$ is general.
\end{lem}

 \begin{proof}
The fact  that $N'$ is perfect follows easily from the fact that $N=N_{C/\P^n}$ is balanced,
with upper subbundle of rank $\equiv k\mod n-1$. For quasi-mixing
 we proceed along the lines of the proof of Lemma 31, Case 2 in \cite{elliptic}. Consider the case
 $n=2m+1$ odd and let $P=P_1\cup _EP_2$ be a fang as there with $E=A_1\times A_2\simeq\P^m\times\P^m$. The factor-flipping gluing of $E=A_1\times A_2\subset P_1$
 to $E=A_2\times A_1\subset P_2$ is obtained by an arbitrary pair of isomorphisms
 $\phi_1:A_1\simeq A_2, \phi_1:A_2\simeq A_1$.\par 
 Start with  the  case $e=2n-1=4m+1$, $k=0$. Consider a reducible curve $C=C_1\cup_xC_2$ where $C_i\subset P_i$
 are birational transforms of  general rational curves $C'_i$ 
 of degree $e_i=n$, such that $C_1\cap C_2=x\in E$, and $p_i\in C_i$ is general.
 Then we have
 \[N_{C_i/P_i}\simeq V^1_i\O(2m+3)\oplus V^2_i\O(2m+2)\]
 where $V^1_i, V^2_i$ are the respective horizontal and vertical subspaces of $T_xE$, i.e.
 $V^1_1=T_xA_1, V^2_1=T_xA_2$ and vice versa (i.e. horizontal for $P_1$ is vertical for $P_2$ and vice versa).
%
  Then we can write
 \[N_{C/P}=R^1\oplus R^2\]
 where
 \[R^1|_{C_1}=V^1_1\otimes\O(2m+3), R^1|_{C_2}=V^1_2\otimes\O(2m+2)\]
 and similarly for $R^2$, where the isomorphism $V^1_1\simto V^1_2$ that
 yields the gluing at $x$ is induced by $\phi_1$, so is a general map up to scalars.
 { 
 Now while $R^1|_{C_1}$ is determined as the upper subsheaf of $N_{C_1/P_1}$, 
 as is its fibre $W_{1,p_1}=R^1|_{p_1}$,
 $R^1|_{C_2}$ on the other hand is only determined at $x$ and its fibre at a general point $p_2\in C_2$
 can be any given general subspace $W_{1,p_2}$ (cf. \S \ref{balance-p1-sec}). 
 Ditto for $W_{2, p_1}, W_{2, p_2} $. Thus, $\tau_{p_1,p_2}$
 is general on $W_{1,p_1}$ while $\tau_{p_1,p_2}\inv$ is general on $W_{2,p_2}$.
 }
  Hence we get that the transport map
  $\tau_{p_1, p_2}$ is quasi-general in this base case $e=2n-1$. \par
 
{\color{red} /************
**** 
}
 ******/\par
 For general degree $e\geq 2n-1$ we use a similar degeneration taking $e_1=e-n+1, e_2=n$ and the
 modifications yielding $N'$ occurring on $C_1$.
 Note that $k$ is congruent to  the remainder of $(n+1)e_1-2$ $ \mod n-1$ so we again assume
 that $N'|_{C_1}$ is of the form
 $m\O(a + 1) \oplus m\O(a)$ for suitable and the argument proceeds as above,
 using the generality of the gluing along $E$. This completes the proof of (i)  for $n=2m+1$ odd
 and $e\geq 2n-1$.
 \par
 {\color{blue}
 Next for $e\geq 4n-3$ we use the same degeneration with $C_2$ of degree $e_2=2n-1$
 and with $C_1$ of degree $e_1=e-2n+2\geq 2n-1$. Then by Case (i) we have that
 $\tau_{p,x}$ is quasi-general, as is $\tau_{x, q}$, hence their composition $\tau_{p,q}$ 
 is general. This proves (ii) for $n$ odd.
 \par
 Finally we consider the case $n=2m$ even.  Here we use an similar degeneration albeit with   
 $A_1=\P^m, A_2=\P^{m-1}$. 
We take $C_1, C_2$ of degrees $e_1=e-2m+1, e_2=2m$, so that $k\equiv 2e_1-2\mod 2m-1$, 
and put the modification centers on $C_1$. Then we have, for suitable $a$,
\[N'_{C_1/P_1}=m\O(a)\oplus (m-1)\O(a-1), N_{C_2/P_2}=(m-1)\O(2m+2)\oplus m\O(2m+1)\]

}
{\color{red} 
 }\par
 {\color{blue}
  As above, 
 at $x$ the upper subbundles will be complementary and glue to horizontal
  subspaces on the other side which extend to
 lower subbundles with general fibre at points $p_1, p_2$.
 }
 \end{proof}
 The result extends to the balanced case:
 \begin{prop}\label{transport-prop}
  Let $p\neq q$ be points in $\P^n, n\geq 3$ and
 let $C$ be a general rational curve through $p,q$  of degree $e\geq 2n-1$
 Let $N'$ be a general simple
 down modification of $N_{C/\P^n}$.  
 Then $N'$ is balanced and quasi mixing. If $e\geq 4n-3$, $N'$ is mixing.
 \end{prop}
 \begin{proof}
 {\color{blue}
 Let $k$ be the remainder of $\deg(N')\mod n-1$. When $k=0$ we are reduced to 
 Lemma \ref{transport-Pn} above. In the general case, it firstly follows as above that $N'$
 is balanced with upper subbundle $N'_+$ of rank $k$. Now let $N"\subset N'$ be a general 
 simple down modification of colength $k$. Then $N"$ is perfect, hence quasi-mixing. 
 We can write
 \[N'=\bigoplus \limits_1^k L_i\oplus\bigoplus\limits_1^{n-1-k}M_j, \]
 with \[\deg(L_i)=a+1, \deg(M_j)=a,
 \]\[N^"=\bigoplus \limits_1^k L_i(-p_i)\oplus\bigoplus\limits_1^{n-1-k}M_j
 \]
 Moving by a suitable projective transformation fixing $p,q$, we may assume $N'_+|_p\subset N'_p$
 is a general subspace, and hence that the
 \[v_i:=L_i|_p,i=1,..k;  w_j=M_j|_p, j=1,...,n-1-k\] form a general basis mod scalars for $N'|_p$ which may naturally
 be identified with a general basis mod scalars for $N"|p$. Then by the perfect case, 
 by varying $(C, N')$ we get $L_i(-p_i)|_q, M_j|_q$ are a general basis mod scalars for $N"|_q$,
 hence $L_i|_q, M_j|_q$ are a general basis mod scalars for $N'_q$. Thus, $\tau_{p,q}$ is quasi general
 on $N'_+|_p$ and $(N'/N'_+)|_p$ so $N'$ is quasi-mixing.\par
 The argument for mixing if $e\geq 4n-3$ is as in the proof of Lemma \ref{transport-Pn}.
 }{\color{red}***


}
 \end{proof}
  We now establish a similar quasi mixing result for  the restricted tangent bundle:
   \begin{lem} \label{transport-ambient}
   Let $p\neq q$ be points in $\P^n, n\geq 3$ and
  let $C$ be a general rational curve of degree $e=kn, k\geq 1$  
  through $p,q$, with restricted tangent bundle $T=T_{\P^n}|_C$.  
  Then $T\simeq n\O_C(k(n+1))$ is perfect and has quasi general parallel transport
   \[\tau_{p,q}:T|_p\to T|_q.\]
   If $k\geq 2$, $\tau_{p,q}$ is general.
    \end{lem}
  \begin{proof} Consider first the case of a rational normal curve $C$ through $p,q$, i.e. $k=1$.
  The family of such curves is acted on transitively by the group $G$ of
  projective transformations fixing $p,q$ and $G$ acts compatibly on the tangent
  spaces to $\P^n$, $T_p, T_q$

 ************** 
  *******\par
  
We  consider a specialization of a rational normal curve in the form of  a polygon $C$ of degree $e=n$, 
i.e. a general chain of $n$ lines
\[C=M_1\cup_{p_1}M_2\cup_{p_2}...\cup_{p_{n-1}} M_n\] with $p=p_0\in M_1, q\in M_n$.
The family of these is an orbit under the action of $G$. 
Then it is easy to see that we can write
\[T=\bigoplus\limits_{i=1}^{n} L_i,\]
\[ L_i=\O_{M_1}(1)\cup...\cup\O_{M_i}(2)\cup...\cup_{}\O_{M_n}(1) , i=1,...,n.\]
where $L_i|_{M_i}=T_{M_i}$ has degree 2 and $L_i|_{M_j}, j\neq i$
has degree 1. 
The above decomposition is
a specialization of a decomposition of the restricted tangent bundle on a smooth
rational normal curve as $n\O(n+1)$ while the associated parallel transport map
specializes to one sending each $L_i|_p$ to $L_i|_q$.
Consider $L_1$ to begin with. Then 
\[L_1|_p=T_pM_1,L_1|_{p_1}=T_{p_1}M_1,
L_1|_{p_2}=T_{p_2}\overline{p, p_2},...,L_1|_{q}=T_q\overline{p,q}.\] Thus,
$\tau_{q,p}=\tau_{p,q}\inv$ sends the fixed 1-dimensional subspace $T_q\overline{p,q}$ to the general subspace
$T_pM_1$.\par
Similarly analyzing $L_2$,  we have
\[L_2|_{p_1}=T_{p_1}M_2, L_2|_p=T_p\overline{p,p_2},...,
L_2|_q=T_q\overline{q, p_2}.
\]
***\par
Thus $\tau_{p,q}$ sends the subspace $T_p\overline{p,p_2}$
to the general subspace $T_q\overline{p_1, q}$, etc.
Hence, moving the polygon fixing $p,q$, we get a quasi-general parallel
transport $\tau_{p,q}$. Therefore the same holds for a general smoothing of $C$ through $p,q$,
which is a general rational normal curve through $p,q$.\par
Then the general case with $e=kn$ follows by induction on $k$ or by considering
a general chain of $k$ rational normal curves 
\[C=C_1\cup_{p_1}C_2\cup...\cup_{p_{n-1}}C_n,\]
with $p=p_0\in C_1, q=p_n\in C_n$ general. Because the parallel transports
 $\tau_{C_i, T, p_{i-1}, p_i}$ are quasi-
general for $i=1,...,n$, so is their composition $\tau_{C, T, p, q}$, and it is even general for $k\geq 2$
by an argument as in the proof  of Lemma  \ref{transport-Pn}.
  \end{proof}
\subsection{Applications}
  The usefulness of having quasi-general parallel transport and holonomy begins to reveal itself in the
    following result (more substantial applications will be given later). 
    We will use the following terminology. If $L$ is a line
    bundle on a reducible curve $C_0=C_1\cup C_2$, its \emph{bidegree} is
    by definition the pair $(\deg(L_{C_1}), \deg(L_{C_2}))$.
    In case $C_0=C_1\cup_{p,q}C_2$ is a  fish curve, any line
    bundle $L$ on $C$ is given by a pair of line bundles $L_i$ on $C_i, i=1,2$
    together with isomorphisms $L_1|_p\simeq L_2|_p, L_1|_q\simeq L_2|_q$.
  \begin{cor}\label{semisimple}
  Let $C_1, C_2$ be curves as in Lemma \ref{transport-ambient} 
  of degrees $k_1n, k_2n, k_1, k_2\geq 1$, both containing the same  $p,q$, 
  and let $C_0=C_1\cup C_2$.
  Then
  \[T_{\P^n}|_{C_0}=\bigoplus L_i\]
  where $(L_1, ...,L_n)$ is a general $n$-tuple of line bundles on $C_0$,
  of bidegree $((n+1)k_1, (n+1)k_2)$.\par
  In particular, if $C$ is a general elliptic curve on degree $kn, k\geq 2$ in $\P^n$ then
  $T_{\P^n}|_C$ is semistable.
  \end{cor}
  \begin{proof}
  By Lemma \ref{transport-ambient}, the holonomy map
  \[\tau=\tau_{C_1, p,q}\inv\tau_{C_2, p.q}:T|_p\to T|_p\]
  is a composite of quasi general maps, hence general, hence diagonalizable with distinct, general  eigenvalues $\alpha_1,...,\alpha_n$.
  Then $T|_{C_0}\simeq \bigoplus L_i$ where $L_i$ is obtained
   by gluing $\O_{C_1}(k(n+1))$ to $\O_{C_2}(k(n+1))$ by the factor $\alpha_i$.
  \end{proof}
  The curve $C_0$ above may be called a perfect fish curve. 
  Smoothing out  to a smooth elliptic curve, the latter result has the following
  consequence: 
  \begin{cor}\label{elliptic-finite-cor}
  If $C\subset\P^n$ is a general elliptic curve of degree $kn, k\geq 2$
  then the set of line subbundles of $\wedge^iT_{\P^n}|_C$ of degree $ik(n+1)$
  is finite.
  \end{cor} 
  \begin{proof}
  Note that for $C_0$ as above,
  \[\wedge^iT_{\P^n}|_{C_0}=\bigoplus\limits_{(j_1<...<j_i)} L_{j_1}\otimes...\otimes L_{j_i}\]
  which is a direct sum of general, distinct line bundles.\par
  Considering a general 1-parameter smoothing $\cC$ of $C_0$ as above, and a line subbundle
  $A_\eta$ of the restriction of $T^i=\wedge^i T_{\P^n}$ on $A_\eta$.
  Arguing as in the proof of Lemma \ref{semistable-lem} we may, after base-changing and blowing up
  that $\cC$ is a smooth fibred surface with special fibre 
  $C_0'=C_1\cup C_2\cup \tilde P\cup\tilde Q$. 
  where $\tilde P, \tilde Q$ are rational chains contracting tp $p,q$ and $A_\eta$
  extends to a line subbundle $A\subset \wedge^iT_{\P^n}|_{\cC}$. 
  As $\wedge^iT_{\P^n}|_{\tilde P\coprod\tilde Q}$ is a trivial bundle, 
  $A$ must have nonpositive  degree on every component of  $\tilde P, \tilde Q$. Therefore clearly
  \[\deg(A_{\tilde P})=\deg(A_{\tilde Q})=0,\]
  so $A$ is trivial on $\tilde P\coprod\tilde Q$, $A_{C_0'}$ descends to $A_{C_0}$ and
  \[\deg A|_{C_i}=ik_j(n+1), j=1,2\]
  so $A_{C_0}$ is one of the $L_{j_1}\otimes...\otimes L_{j_i}$ 
  with the $L_i$ as above. By semicontinuity, the Corollary follows.
  
 \end{proof}
 This result implies semistability of the restricted tangent bundle for perfect elliptic curves; 
 in the next section we extend this result to all elliptic curves of large degree (see
 Theorem \ref{elliptic-ambient-semistable}). 
  \par
 As a consequence of Lemma \ref{transport-Pn}, 
 we deduce an analogous    result about parallel transport for rational curves in a fang component.
 This situation is slightly different in that the transfer is between points are not general on the curve.
 \begin{lem}\label{curve-in-fang-lem}
 Let $P=B_{A_1}\P^{n}, A_1=\P^k$ with exceptional divisor $E=A_1\times A_2$,
 let $p,q\in E$ be general and let $C\subset P$
 be a general curve with $C.E=\{p,q\}$ that is
 the birational transform of a rational curve $C'\subset\P^{n}$ of degree $e\equiv 1\mod 2m, e\geq 4m+1$.
 Let $N=N_{C/P}$ and $\tau=\tau_{C, N, p,q}$. \par
 (i) Assume $n=2m+1,  k=m, e\equiv 1\mod 2m$.
 Then  the parallel transport $\tau$ is quasi general.
 \par (ii) Assume $n=2m, k=m-1, e\equiv 1\mod 2m-1, e\geq 4m-1$. Then the parallel transport on the
 upper subbundle $N^+\subset N$, of rank $2m-2$,  is quasi-general.
 \end{lem}
 NB: $p,q $ are evidently not general on $C$.
 \begin{proof} We will prove (i) as the proof of (ii) is similar.
 {\color{blue}
 Firstly, $N'=N_{C'/\P^n}$ is perfect and 
 the corresponding points $p',q'$ are general on $C'$ so by Lemma \ref{transport-Pn}
 $\tau':=\tau_{N', p', q'}$ is quasi-general. Now $N'$ is a down modification of $N$ at $p,q$
 and we can decompose
 }
 {\color{blue} 
 \[N=N_1\oplus N_2=(V_p\oplus V_q)\otimes\O(e-1+(e-1)/m)\]
where $V_p=T_pA_1, V_q=T_qA_1$ and the map $N\to N_{C'/\P^n}=2m\O(e+(e-1)/m)$
when restricted on $N_1$ is fibrewise injective at $p$ and zero at $q$, and vice versa for $N_2$,
i.e. $N_1$ is horizontal at $p$ and vertical at $q$ and vice versa for $N_2$.
 }
 Now as remarked, $\tau'$ is quasi general
 hence quasi-general on $V_p$, while $(\tau')\inv$ is general on $V_q$. 
 Now coordinates on $P$ have the form $x_i/x_j$ with $x_i, x_j$ coordinates on
 $\P^{2m+1}$ and $x_j$ vanishes on $A_1$. Then  $\tau $ is quasi-general on $V_p$
 while $\tau\inv$ is quasi general on $V_q$.  Hence $\tau $ is quasi-general.
 
 \end{proof}
 \begin{lem}\label{qg4down-lem}
 
 Notations as in Lemma \ref{curve-in-fang-lem} (no congruence assumptions on $e, n$), 
 let $N'\subset N$ be a general simple
 modification of colength equal to the remainder of $\deg(N) \mod n-1$. Then
 $N'$ is perfect and $\tau_{C, N', p,q}$ is quasi-general.
 \end{lem}
\begin{proof}
For a suitable general splitting into line bundles $N=\bigoplus L_i$, we can write
\[N'=\bigoplus\limits_1^k L_i(-p_i)\oplus\bigoplus\limits_{k+1}^r L_i.\]
Consequently $\tau_{C, N', p,q}$ can be identified with $\tau_{C, N, p,q}$, so
it is quasi-general.
\end{proof}
 \par
{\color{blue}The following result related the 'generality' of a subbundle to its degree:
will be useful in the proof of Theorem \ref{g=1-thm}.
\begin{lem}\label{subbundle-lem}
On $\P^1$, let $F\subset  E=s\O(d)$ be a subbundle of rank $r$ such that 
\eqspl{dimension-eq}{\dim (\tau_{E, p,q}(F|_p)\cap F|_q)\leq r-a.}
Then $\deg(F)\leq rd-a$.
\end{lem}
\begin{proof}
We can write $F=\O(b_1)\oplus...\oplus\O(b_r), d\geq b_1\geq...\geq b_r$. If $b_{r-a+1}=d$
it would violate \eqref{dimension-eq}. Therefore $d-1\geq b_{r-a+1}\geq...\geq b_r$
so $\deg(F)\leq dr-a$. 
\end{proof}
}
	\section{Genus 1 in Projective space}
	As mentioned above,  Ein-Lazarsfeld \cite{ein-laz-normal}
	have proven c-semistability of the normal bundle of an elliptic normal curve, of degree $n+1$ in $\P^n$. 
	Here we prove an analogous but logically independent result, showing semistability of the normal
	bundle of a general elliptic curve of degree  $\geq 3n-3$ in $\P^n, n\geq 4$.
In the case $n=3$, 	the analogous result  follows from 
		Coskun-Larson-Vogt \cite{coskun-larson-vogt-normal}.\par
	First, as a matter of terminology, the \emph{bidegree} of a bundle $E$
	on a reducible curve $C_1\cup C_2$ is by definition $(\deg(E|_{C_1}), \deg(E|_{C_2}))$.
	A down or up modification is \emph{simple} if it is locally of corank 1.
	\begin{thm}\label{g=1-thm} The normal bundle $N$ of
	a general elliptic curve of degree  $e\geq 4n-4$ in $\P^n, n\geq 3$ is semistable as is
	any general simple modification $N'\subset N$.
	Moreover if the slope  $\mu(N')$ is not
	an integer then $N'$ is stable.
		\end{thm}
	In view of Lemma \ref{semistable-lem}, this implies:
	\begin{cor}
	Notations as above, any smooth $m$-dimensional subvariety $Y$ containing $C$ must have
	\[C.(-K_Y)\leq\frac{(n+1)(m-1)e}{n-1}.\]
	\end{cor}
To prove Theorem \ref{g=1-thm}, we begin with a useful lemma. Note that if $E=L_1\oplus L_2$ with $L_1, L_2$ \emphred{distinct} line bundles 
of degree 0 and  $E'\subset E$
is a \emphred{general}  modification of colength 1 then any line subbundle $M\subset E'$
must have degree $\leq -1$ hence $E'$ is stable. This is false for non-general modifications
or if the line bundles are not distinct. 
The following Lemma extends
this to any rank.
\begin{lem}\label{downmod-lem}
Let $C$ be a nodal curve of arithmetic genus 1, and let $E=\bigoplus\limits_1^rL_i$
be a direct sum of a general collection of  line bundles of degree 0. \par
 (i) Let $E\to E^+\to\bigoplus\limits_1^k \tau_{p_j}$
be a general (independently of the $L_i$) simple up modification with  $0<k<r$. Then
$E^+$ is stable.\par
(ii) Let $E'\to E\to\sum\limits_1^k  \tau_{p_i}$ be a general simple down modification, $0<k<r$.
Then 
 $E'$ is stable.
\end{lem}
\begin{proof}[Proof of Lemma] By duality it suffices to prove (ii).
{\color{red}/**********
%
%
%
%
%
*********/
}
\par 
Let $F'\subset E'$ be a subbundle of rank $s<r$ and let $F\subset E$ be its saturation
and set $-a=\deg(F)$. As discussed in \S \ref{quot-sec}, deformations of $F$ as subbundle of $E$
are parametrized by the appropriate Quot scheme \cite{grothendieck-hilbert}, \cite{sernesi}. 
Note that if $F$ contains any summand $L_i$ then $L_i\cap E'=L_i$,
so by modding out by this $L_i$ and using  induction on $r$, we may assume this does not happen.
Then  we have exact
\[\exseq{F'}{F}{\tau_0}\]
where by rearranging notation we may assume $\tau_0=\sum\limits_1^b\tau_{p_i}$
where each $\tau_{p_i}=k(p_i)$ has length 1. 
Then we have
\[\deg(F')=-a-b.\]
If $b=0$ then $F'=F$ and clearly $\mu(F)<\mu(E')$, so assme $b>0$.
Let $h_i\subset E|_{p_i}$
be the hyperplane defining the torsion quotient $\tau_{p_i}$. Then we have
\[F|_{p_i}\subset h_i, i=b+1,...,k.\]
This is $s$ conditions on $F$ for each $i$, so $(k-b)s$ conditions in total
(a bit more formally, we may in a standard way construct the space of tuples $(F, p_\bullet, h_\bullet)$ 
as fibred  over the space of $(p_\bullet, h_\bullet)$ and consider its general fibre).
By induction on $r$ as noted, we may assume $F$ does not contain any $L_i$ summand, i.e. 
$h^0(F(-L_i))=0, \forall i$, hence
$h^1(\check{F}(L_i))=0, \forall i$ so that $h^1(\check F\otimes E)=0$ and consequently
\[h^1(\hom(F, E/F))=0.\]
Now in general, if $U, V$ are vector bundles of respective ranks $u,v$ then
\[c_1(U\otimes V)=uc_1(V)+vc_1(U).\]
In our case this shows that $\deg(\hom(F, E/F))=ra$. So we conclude
that the family of bundles $F$ is unobstructed of dimension
\[\dim(\{F\})=\chi(\hom(F, E/F))=ra.\] 
Because the points $p_i$ and hyperplanes $h_i$ are general,
i,e, we are considering a general fibre of $\{(F, p_\bullet, h_\bullet)\}\to \{(p_\bullet, h_\bullet)\}$,
it follows that
\[(k-b)s\leq ar.\]
Consequently
\[\mu(F')=\frac{-a-b}{s}< \frac{-ar-bs}{rs}\leq\frac{-k}{r}=\mu(E'),\]
showing stability.

{\color{red}/**********
%
%
%
}
\end{proof}
	\begin{proof}[Proof of Theorem] Write $\deg(N')=\alpha(n-1)+a, a\in [0,n-2]$.
	Replacing $N'$ by a twist, we may assume $\alpha\geq 4n$. 
	Assume first that $n=2m+1$ is odd. \par
		Consider a fang degeneration
		\[P_0=P_1\cup_QP_2\]
		where
		\[P_i=B_{\P^m}\P^n\supset Q_i=\P^m\times\P^{m}\]
		in which $Q_i\subset P_i, i=1,2$ are exceptional divisors and 
		$P_0$ is constructed via isomorphisms $Q_1\simeq Q\simeq Q_2$ which may be assume general. 
		There is a standard smoothing of $ P_0$ to $\P^n$.
		Consider
		curves
		\[C_1\subset P_1, C_2\subset P_2\]
		with each being a birational transform of a general rational curve
		of degree $e_1$ resp. $e_2$ meeting $\P^m$ resp. twice, such that
		\[C_1\cap Q=C_2\cap Q=\{p,q\}.\] 
		Specifically we take $e_1=4m-2, e_2\geq 4m-1$, so that $N_{C_1/P_1}$ is perfect while
		$N_{C_2/P_2}$ is balanced.  Let $N_2$ be a general modification of $N_{C_2/P_2}$
		and $N'=N_1\cup N_2$ (i.e. a modification of $N$ centered on $C_2$) .Thus
		\[N_1:=N_{C_1/P_1}=2m\O(d_1), \]
		\[N_2:=a\O(d_2+1)\oplus (2m-a)\O(d_2), 0\leq a\leq m\]
		so $\deg (N')=2m(d_1+d_2)+a$ and $d_1+d_2=[\mu(N')]$.
		{\color{blue}
		Let $N_2^0\subset N_2$ be a general down modification of colength $a$.
		Then $N_2^0$ is perfect and by Lemma \ref{qg4down-lem}, $\tau_{C_2, N^0_2, p,q}$ is quasi-general.
		Consequently, $\eta= \tau_{C_1, N_1, p,q}\circ\tau_{C_2, N_2^0, q, p}$ is general and
		we can decompose
		\[N^0:=N_1\cup N_2^0=\bigoplus L_i\]
		with the $L_i$ a general collection of line bundles of the same degree, corresponding to 
		eigenvectors of the holonomy $\eta$. Now  $N'=N_1\cup N_2$ is an up modification of $N^0$, and by varying $C_1$
		we can ensure that this modification is general relative to the eigenvectors, i.e. to the
		above decomposition of $N^0$.  Therefore by Lemma \ref{downmod-lem}, $N'$ is semistable.
		
		}
		{\color{red}
		/***********
%
%
***********/	}
	\end{proof}

	Next we prove an analogous result for the restricted tangent bundle:
	 \begin{thm}\label{elliptic-ambient-semistable}
	  If $C\subset\P^n$ is a general elliptic curve of degree $e\geq 2n$ then $T_{\P^n}|_C$ is
	  semistable.
	  \end{thm}
	  \begin{proof}
	  We will use the setup and method of Corollaries \ref{semisimple} and \ref{elliptic-finite-cor}
	  and Lemma \ref{downmod-lem}.
	  Consider a reducible curve in the form of a hairy fish:
	  \eqspl{c*}{C_*=C_0\cup M_1\cup...\cup M_m \subset\P^n}
	  with $C_0=C_1\cup _{p,q}C_2$ a fish curve of arithmetic genus 1 and total degree $kn, k\geq 2$, 
$m$ is the remainder of $e$ mod $n$	  and the $M_i$ 
are general 1-secant lines of $C_0$ for $i=1,...,m$.
	  As proven above in 
	  Corollary \ref{semisimple}, we have \[T_0=T_{\P^n}|_{C_0}=\bigoplus L_j\]
	  for a general $n$-tuple of line bundles $L_1,...,L_n$ of total degree $k(n+1)$ each, while
	  {\color{blue}\[T_j=T_{\P^n}|_{M_j}=\O(2)\oplus(n-1)\O(1)\]
	  where the  $\O(2)$ summand corresponds to $T_{M_j}$, hence
	  its fibre at the node $p_j=M_j. C_0$ is in general position with respect to
	  those of the $L_i$. Let 
	  \[F\subset T_*:=T|_{C_*}\] be a subbundle of rank $s$ and
	  assume $F$ does not contain $L_i$ for any $i$ as the contrary case is
	  similar and simpler. Let 
	  $F_0=F|_{C_0}, a=\deg(F_0)$, and suppose that $F|_{M_j}$ contains $\O_{M_j}(2)$
	  for $j=1,...,\l$ only. Then
	  \[\mu(F)\leq 1+(a+\l)/s.\]
	  On the other hand, as in the proof of Lemma  \ref{downmod-lem}, $F_0$ moves in an
unobstructed family (Quot scheme, cf. \S \ref{quot-sec}) of subbundles $F_0\subset T_0$ of
	  dimension 
	  \[\dim(\{F_0\})=\chi(\hom(F_0, T_0/F_0))=skn(n+1)-a\] 
	  and it is $n-s$ conditions for $F_0$ to contain
	  each $T_{p_j}M_j$, hence $\l(n-s)$ conditions in total,  by dimension-counting as
	  in the proof of Lemma \ref{downmod-lem}, we have
	  \[\l(n-s)\leq skn(n+1)-an.\]
	  Therefore
	  \[\mu(F)\leq 1+(a+\l)/s\leq 1+\l/n+k(n+1)\leq 1+m/n+k(n+1)=\mu(T_*)\]
	   hence $T_*$
	  	   is semistable.
	  }
	  
	  {\color{red}***********}

	  \end{proof}
	  The following extension of Theorem \ref{elliptic-ambient-semistable} will be useful in the sequel:
	  \begin{lem}\label{modified-tangent-lem}
	  Let $C\subset\P^n$ be a general elliptic curve of degree $e\geq 2n$ and let
	  $T'\subset T=T|_{\P^n}|_C$ be a modification locally of corank $\leq 1$
	  that fits in an exact sequence
	  \[\exseq{T'}{T}{A}\]
	  with the properties\par
	  (i) the support  of $A$ is contained in a general hypersurface section $F_d.C$;\par
	  (ii) the map $T\to A$ is general;\par
(iii)  $A$ can be decomposed as $A=A_1\coprod A_2$ where
	  the total  length of $A_1$ is divisible by $n$ (possibly 0) and the support of $A_2$ consists of
	  at most $d(e-n)$ points.\par
	  Then $T'$ is semistable.
	  \end{lem}
	  \begin{rem}
	  What's different  here is that the support of $A$ is not necessarily a general set of points.
	  \end{rem}
	  \begin{proof}
	  Consider a degenerate form $C_0=C_1\cup_{p,q}C_2$ with $C_1, C_2$ rational
	  curves, $C_1\subset\P^n$ a rational normal curve, $A_1\subset C_1, A_2\subset C_2$,
	   and $\tau_{T, C_1, p,q}$ general.
	  Then $T'|_{C_1}$ is perfect while $T'|_{C_2}$ is balanced, so the above argument applies.
	  \end{proof} 
	\begin{example}
	The case we require in the sequel (proof of Theorem \ref{hypersurface-ambient-thm}) is when
	 $A=kF_d.C$, i.e. $A$ has the same multiplicity $k$
	at every point of $F_d.C$.
	\end{example}

	\section{Higher genus in Projective space}
	The following result extends Theorem \ref{g=1-thm} to higher genus,
	using an induction on the genus. 
	\begin{thm}\label{genus g}
	Let  $C$ be a general curve of genus $g\geq 2$ and degree $e\geq 4g(n-1)-g+1$
	in $\P^{n}$ and let $N$ be its normal bundle. Then
   $N$ is hyper-stable;\par

	 	\end{thm}
	\begin{cor}
	Notations as above, if $Y\subset\P^n$ is a smooth $p$-dimensional variety containing $C$ then
	\[C.(-K_Y)\leq((p-1)(n+1)e-(n-p)(2g-2))/(n-1).\]
	\end{cor}
	\begin{proof} [Proof of Theorem] The proof is by induction on $g$, based on the case $g=1$ 
	proved above (Theorem \ref{g=1-thm}).
	We consider a fang degeneration as above with
	\[P_0=P_1\cup_E P_2, E=\P^m\times\P^{n-1-m},\]
	with $m$ to be determined, and a lci curve of the form
	\[C_0=C_1\cup_p C_2\]
	with $C_1\subset P_1$ the birational transform of a general curve of 
	genus $g-1\geq 1$ and degree $e_1\geq4(g-1)-(g-1)+1$, and $C_2\subset P_2$ the birational
	transform of a general curve $C_2'$ of genus 1,and degree $e_2\geq 4n-1$,
	such that $2e_2-2-m\not\equiv 0\mod n-1$. The latter incongruence ensures that
	the slope of $N_{C_2/P_2}$, which is a general down modification of colength $m$
	of $N_{C'_2/P_2}$, is not an integer, hence by Theorem   \ref{g=1-thm}  $N_{C_2/P_2}$ is stable.
	Likewise, by induction or by Theorem \ref{g=1-thm}, $N_{C_1/P_1}$ is hyper- semi-stable.

\par
Now the following  elementary Lemma (just a restatement of  Lemma \ref{openness-lem})
applied to a general down modification
of $N_{C_0/P_0}$ centered on $C_1$ shows that  $N_{C_0/P_0}$ is, in a suitable sense
hyper-stable, hence so is its smoothing is a smoothing $N_{C/P}$.
	\begin{lem}\label{union-lem}
	Let $E_0$ be a vector bundle on a connected nodal curve $C_0$ that is the union
	of smooth components, such that the restriction of $E_0$ on each component is
	semistable (resp. c-semistable). Then for any smoothing $(C, E)$ of $(C_0, E_0)$, we have\par 
	(i)  $E$ is semistable (resp. c-semistable);\par
	(ii) if moreover the restriction of $E_0$ on at least one component of $C_0$ is stable (resp. c-stable), 
	$E$ is stable (resp. c-stable).
	\end{lem}\end{proof}
	Next we prove the higher-genus case of Theorem \ref{elliptic-ambient-semistable}.
	\begin{thm}\label{ambient-semistable-thm}
	If $C\subset \P^n$ is a general curve of genus $g\geq 1$ and degree $e\geq 2gn$
	then\par
	(i)  the restricted tangent bundle $T=T_{\P^n}|_C$ is c-semistable;\par
	(ii) if $T'\subset T$ is a modification as in Lemma \ref{modified-tangent-lem}, then $T'$ is
	semistable.
	\end{thm}
	\begin{proof}
	(i) follows easily from Theorem \ref{elliptic-ambient-semistable} and Lemma \ref{union-lem}
	using either induction or a chain of elliptic curves. (ii) follows similarly using Lemma \ref{modified-tangent-lem}.
	\end{proof}

	\section{Fano hypersurfaces}
	The purpose of this section is to construct curves with (semi)stable normal
	or restricted tangent  bundle
	on some general hypersurfaces of dimension $n\geq 3$ and degree $d\leq n+1$ 
	in projective space.
	We begin with the case of anticanonical hypersurfaces (degree $n+1$ in $\P^{n+1}$),
	where we are able to construct such curves of all sufficiently large degrees.
	\begin{thm}\label{hypersurface-thm}
Let $X$ be a general hypersurface of degree $n+1$ in $\P^{n+1}, n\geq 4$ and 
let $ g, e$ be such that $e\geq n( 4g(n-1)-g+1)$. 
	Then  if $g\geq 2$ (resp. $g=1$), $X$ contains a curve of genus $g$ and degree $e$ with hyper- stable 
	(resp. hyper semi-stable) normal bundle in $X$.
	\end{thm}
	\begin{rem}
	Note that for $C, X$ as above the slope 
	\[\mu(N_{C/X})=\frac{e(n+1)+2g-2}{n-1}\equiv \frac{2(e+g-1)}{n-1}\mod \Z \] is generally not an integer.
	\end{rem}
	\begin{cor}
	Notations as above, if the curve $C$ is contained in a smooth $p$-dimensional subvariety $Y\subset X$
	then
	\[C.(-K_Y)\leq((p-1)(n+1)e-(n-p)(2g-2))/(n-1).\]
	\end{cor}
	\begin{proof}[Proof of Theorem]
	We will use the same fan- quasi cone degeneration as in \cite{caudatenormal}. Thus we take
	\[P_0=P_1\cup_Q P_2\]
	with $P_1$ the blowup pn $\P^{n+1}$ at a point $p$, with exceptional divisor $Q=\P^n$,
	and $P_2=\P^{n+1}$ containing $Q$ as a hyperplane. In $P_0$ we consider a Cartier divisor
	\[X_0=X_1\cup_Z X_2\]
	with $X_1$ the blowup of a quasi-cone $\bar X_1$, i.e.
	a hypersurface of degree $n+1$ with multiplicity $n$ at $p$,
	and $X_2$ a general hypersurface of degree $n$ in $P_2=\P^{n+1}$. Then
	$X_0\subset P_0$ smooths out to a hypersurface $X\subset\P^{n+1}$ of degree $n+1$.
	As noted in
	\cite{caudatenormal}, \S 4, $X_1$ may be realized as the blowup of $\P^n$ in an $(n,n+1)$
	complete intersection 
	\[Y=F_n\cap F_{n+1}\] where $F_n$ corresponds
	to $Z=X_1\cap Q=X_2\cap Q$ and $F_{n+1}$ (which is not unique)
	corresponds to a hyperplane section of  $\bar X_1$. 
	 Proceeding as in \cite{caudatenormal}, write
\[e=(n+1)e_1-a, 0\leq a\leq n,\]	 and consider a curve $C_0\subset X_0$
of the form 
\[C_0=C_1\cup C_2\]
where the components are as follows.\par
 $C_1\subset X_1$ is the birational (=isomorphic) transform of a curve $C'$ of genus $g$
and degree $e_1$ in $\P^n$, meeting $Y$ in $a$ points $p_1,...,p_a$ with general tangents, whose normal bundle
$N_{C'/\P^n}$ is stable (resp. semistable, for $g=1$),
and remains stable (resp. semistable) after the general down modification
at $p_1,...,p_a$, corresponding to the tangent spaces $T_{p_i}Y$. Such a curve exists
by Theorem  \ref{genus g} and meets $Z$ in $C'\cap F_{n}\setminus\{p_1,...,p_a\}$. 
The latter modification coincides with $N_{C_1/P_1} $.\par
As for $C_2\subset X_2$,  as in \cite{caudatenormal}
it is  is a disjoint union of lines with trivial 
(hence semistable) normal bundle, meeting $C_1$ in $C_1\cap Z$. 
Now we have
\[N_{C_0/P_0}|_{C_i}=N_{C_i/P_i}, i=1,2.\]
Therefore e.g. by Lemma \ref{union-lem} above or by an argument as in \cite{caudatenormal},
a smoothing $C\subset X$ of $C_0\subset X_0$ to a curve on a hypersurface
of degree $n+1$ in $\P^{n+1}$ has stable (resp. semistable) normal bundle.\par
{\color{red}
residual meeting in $F_n$ or $F_{n+1}$?
}
	\end{proof}
	Next we prove an analogous semistability results for the restricted tangent bundle.
		\begin{thm}\label{hypersurface-ambient-thm}
	Let $X$ be a general hypersurface of degree $n+1$ in $\P^{n+1}, n\geq 4$. 
	Then for all  $ g\geq 1$ and large enough $e$ depending on $g$ and $n$,
    $X$ contains a curve $C$ of genus $g$ and degree $e$ such that the
    restricted tangent bundle $T_X|_C$ is semistable.
		\end{thm}
		\begin{proof}
		We use that same construction as above. Note (see e.g. \cite{elliptic}, \S 2) that
		a limit on $X_0$ of the tangent bundle of a general hypersurface has
		the form
		\[T_0=T_{X_1}\mlog{Z}\cup T_{X_2}\mlog{Z}.\]
		Now we choose $C'$ so that $C'\cap Y=C'\cap F_n$, hence $C\cap Z=\emptyset$
		and $C_2=\emptyset$.
		Also, $T_0|_{C_1}=T_{X_1}|_{C_1}$. This is the locally corank-1 modification
		of $T_{\P^n}|_{C_1}$ corresponding to the hyperplanes $T_pY+T_pC_1$
		for all $p\in C_1\cap Y$, which we may assume is a general modification.
		Hence we may conclude by Lemma \ref{modified-tangent-lem}.
		\end{proof}
	
	Next we  take up the case of hypersurfaces of degree $d<n$ in $\P^n$, where our results
	are spotty as to the curve degrees:
	\begin{thm}\label{lower}
	Let $X$ be a general hypersurface of degree $d, \ \ (n+1)/2\leq d< n$ in $\P^n$ 
	and suppose that $g\geq 1$,  $0<<e_0\leq e$ are such that
	
	\eqspl{numerical}{
	(d-2)(n-d+1)e-(n-2)de_0+2(n-d)(g-1)=0.
	}

Then  $X$ contains  smooth curve of degree $e$ and genus $g$
	with semistable normal bundle.\end{thm}
	By arithmetics, we conclude:
\begin{cor}	Assume
		\eqspl{divisible}
			{((d-2)(n-d+1), d(n-2)\ \ |\ \ 2(n-d)(g-1).
			}
Then with $X, d, g$ as above, there exists an arithmetic progression of $e$ values for which
$X$  contains a $(g, e)$ curve with semistable normal bundle.\par
	In particular, whenever 
	\eqspl{progression}{d(n-2)\ |\ 2(n-d)(g-1),
	}
	e.g. if $g=1$, 
	such curves exist for all large degrees $e$ divisible by $d(n-2)$.
	\end{cor}
	\begin{rem}\begin{enumerate}\item
	Note that \eqref{numerical} implies that for large $e_0$, we have
	\[e/e_0\sim\frac{n-2}{n-d+1}\frac{d}{d=2}>1.\]
	So $e_0>>0$ implies that $e>e_0$. Conversely $e>>0$ also implies $e_0>>0$.\item
	For $g=1$, \eqref{numerical} reads $(d-2)(n-d+1)e=(n-2)de_0$ which is solvable for
	all sufficiently divisible $e$.\item
	For given $g>1$, \eqref{numerical} admits a solution $(e_0, e)$ iff \eqref{divisible} holds.
Once one solution exists, there will exist solutions with $e$ filling out an arithmetic progression.
	\item In case $d=n-1$, \eqref{divisible} means that either $n$ is even or $g$ is odd.\item
	In case $d=n-2$, \eqref{divisible} means that either $g\equiv 1\mod 3$ or $n\not\equiv 2\mod 3$.

	\end{enumerate}
	\end{rem}
\begin{rem}
Consider a projective bundle $Z=\P(E)$ and a divisor of the form $X=\P(F)\subset Z$ where
$F$ is a torsion-free sheaf that fits in an exact sequence 
\[\exseq{M}{E}{F}\]
with $M$ a line bundle.
Then $X$ is the zero locus of the composite
\[\O_Z(-1)\to \pi^* \check E\to \pi^*\check M\]
i.e.
\[[X]=c_1(\check M\otimes\O_Z(1)).\]
When the zero-locus of the map $M\to E$ has codimension $3$ or more (e.g. when $E$
has rank $\geq 3$ and the map is general and $\check M\times E$ is globally generated), 
then $X$ is irreducible.
\end{rem}
	\begin{proof}[Proof of Theorem]
	We will use the fang setup as in \cite{caudatenormal}, \S 6 and in  \S \ref{fangs} above. 
	Thus we set $m=d-1$
	and consider a limiting form of $\P^n$ which is a fang of the form
	\[Z_0=Z_1\cup Z_2\]
	where
	\[Z_1=\P_{\P^m}(1, 0^{n-m})=B_{\P^{n-m-1}}\P^n, Z_2=\P_{\P^{n-m-1}}(1, 0^{m+1})=B_{\P^m}\P^n.\]
As explained in \S \ref{fangs}, the degeneration may be
	 constructed by blowing up $\P^{n-m-1}\times 0$ in $\P^n\times\A^1$. This yields
	 a smooth variety $\cZ/\A^1$ with general fibre $\P^n$ and special fibre $Z_0$.
	In $Z_0$ we consider a divisor which is a limiting form of a degree-$d$ hypersurface
	in $\P^n$ and has the form
	\[X_0=X_1\cup X_2\]
	where: $X_1=\P_{\P^m}(G)\subset Z_1$ is the blowup of a hypersurface of degree $d$
	in $\P^n$ with multiplicity $m=d-1$ along $\P^{n-m-1}$, which is a (singular)
	fibration over $\P^m$ with  general fibre a hyperplane in a  $\P^{n-m}$ fibre of $\pi_1:Z_1\to \P^m$ 
	that is
	 residual to $(d-1)\P^{n-m-1}$; and $X_2\subset Z_2$ is the blowup of a hypersurface
	 of degree $d=m+1$ containing $\P^m$ and  is fibred over $\P^{n-m-1}$
	with general fibre a general hypersurface of degree $m=d-1$, 
	residual to $\P^m$, in the $\P^{m+1}$ fibre of $Z_2$.
	We recall that $G$  has the form
	\[G=\pi_{1*}(\pi_1^*(\O(1)\otimes (n-m)\O)\otimes\O_{X_1}).\]
	Thus $G$ is a rank-$n-m$ reflexive sheaf  over $\P^m$ which fits in an
	exact sequence
	\eqspl{g}{
	\exseq{\O(-d+1)}{\O(1)\oplus (n-m)\O}{G}.
	} $G$ is locally free except over the locus $W\subset\P^m$ where the left map in
	\eqref{g} vanishes,
	which is a complete intersection of codimension $1+n-m=n-d+2$ (or empty if $m<(n+1)/2$).
	Then in $X_0$ we consider a connected lci curve of the form
	\[C_0=C_1\cup C_2\]
	where  (as in \cite{caudatenormal}, \S 6),
	 $C_2\subset X_2$ is a disjoint union of lines contained in fibres
	of the projection to $\P^{n-m-1}$, with trivial (hence semistable)
	normal bundle in $X_2$, while $C_1\subset X_1$ is a suitable isomorphic  lift of $\P^n$-degree $e$ 
	of a smooth curve $C_+\subset\P^m$
	of genus $g$ disjoint from $W$ with semistable normal bundle.
	Then $C_0\subset X_0$ deforms to a smooth curve $C\subset X$
	of degree $e$ and genus $g$ on a general hypersurface of degree $d$ and the normal bundle 
	$N_{C/X}$ will be semistable if $N_{C_0/X_0}$ is, which in turn will be true provided
	$N_{C_1/X_1}$ is semistable. Thus it would suffice to show that with suitable choices
	$N_{C_1/X_1}$ may be assumed semistable.
	For convenience, let us call the $\P^n$ and $\P^m$ degrees of a curve $C_1\subset X_1$
	the upper and lower degrees, say $e, e_0$, and the pair $(e, e_0)$ the bidegree.
	\par Our strategy for constructing  $C_1$  with semistable normal bundle is based on the following exact sequence
	\eqspl{normal-seq}{\exseq{T_v|_{C_1}}{N_{C_1/X_1}}{N_{C_+/\P^m}}}
	combined with the following easy remark
	\begin{lem}\label{slope-lem}
	Let
	\[\exseq{E_1}{E}{E_2}\]
	be an exact sequence of vector bundles on a curve such that  two of $E_1$ , $E_2$ and $E$ 
	are semistable with the same slope $\mu$. Then  so is the third. 
	\end{lem}  

%
Now given $C_+\subset\P^m$, lifts of $C_+$ to $C_1\subset X_1=\P(G)$
of $\P^n$-degree $e$ correspond to invertible quotients
\[\exseq{K}{G|_{C_+}}{ B}, \deg(B)=e,\]
i.e.  the $\P^n$-degree of $C_1$ coincides with $\deg(B)$.
Then
	 we have  as in loc. cit. $T_v=K^*(B)$ where $K$ is the kernel as above
 and $T_v|_{C_1}$ has slope 
		\[\mu(T_v|_{C_1})=e+(e-de_0)/(n-d).\]
		On the other hand we have \[\mu(N_{C_+/\P^m})=\frac{(m+1)e_0+2g-2}{m-1}
		=\frac{de_0+2g-2}{d-2}=e_0+\frac{2e_0+2g-2}{d-2}. \]
		Now by Theorems \ref{g=1-thm} and \ref{genus g}, we may assume
		$N_{C_+/\P^m}$ is semistable. 
Equating the slopes of $N_{C_+/\P^m}$ and $T_v|_{C_1}$  per Lemma \ref{slope-lem} leads to
		\eqspl{e-equation}{e(d-2)(n-d+1)-(n-2)de_0+(n-d)(2g-2)=0
		}
		The condition for $C_1$ of bidegree $(e, e_0)$ and genus $g$ to exist is that  \eqref{e-equation} 
		should hold with $e_0>>0$ and $e-e_0\geq 0$ and that $T_v|_{C_1}$ is semistable.\par
		We will construct $C_1$ in the form of a comb. To this end,
		we first construct some good curves of genus 0 to serve as shaft.
		Let $C_{10}\subset \P^m$ 
		be a general rational curve of degree $e_0\geq m$. \par
		{\color{red}************
*******************}\par
		By Lemma \ref{g-on-c-lem}, $G|_{C_{10}}$ is balanced (the numerical condition
		there is equivalent to  $d\geq (n+1)/2$). Then choosing a general surjection
		{ ${G|_{C_{10}}\to B}$}
		where $B$ is a line bundle of degree $e_1\geq e_0$, we have by Lemma \ref{kernel-lem} 
		that the kernel $K$
		is balanced, hence for the lifting
		of $C_{10}$ to $X_1$ corresponding to $B$,  
		$T_v|_{C_{10}}$ will be balanced.  In particular if we choose
		$e_1$ so that $(n-d)|(e_1-de_0)$, e.g. $e_1=de_0$, then
		$T_v|_{C_1}$   will be perfect, hence semistable.\par
		Now for given genus $g\geq 1$ consider a connected nodal curve of genus $g$ of the 
		comb-like form
		\[C_1=C_{10}\cup D_1\cup...\cup D_g\]
		Where $C_{10}$, the shaft,  is as above and 
		the $D_j$, the teeth,  are elliptic curves of degree $\epsilon_j\geq 2(n-m)$, each contained
		in a $\P^{n-m}$ fibre of $X_1/\P^m$.
		So $D_j$ has bidegree $(\epsilon_j, 0)$.
		 Thus $C_1$ has $\P^n$ degree
		\[e=e_1+\sum\epsilon_j\geq de_0+\sum\epsilon_j\] 
		and by choosing $\epsilon_j$ suitably  we may arrange things so that $e$ satisfies \eqref{e-equation}. By Theorem \ref{elliptic-ambient-semistable},
		we have that each $T_{\P^{n-m}}|_{D_j}=T_v|_{D_j}$ is semistable, 
		hence so is $T_v|_{C_1}$. 
		Also,  as for $T_{\P^m}$, we have that  $T_{\P^m}|_{C_{10}}$ is semistable 
		while $T_{\P^m}|_{D_j}$ is trivial, hence semistable 
		for each $j$, hence  $T_{\P^m}|_{C_{1}}$ is semistable. 
		Smoothing $C_1$ then proves the Theorem.
		\end{proof}
		\newpage
\bibliographystyle{amsplain}

\bibliography{../mybib}
		
\end{document}